\documentclass[11pt]{amsart}
\usepackage{amsmath}
\usepackage{amsxtra}
\usepackage{amscd}
\usepackage{amsthm}
\usepackage{amsfonts}
\usepackage{amssymb}
\usepackage{eucal}

\newtheorem{cor}[subsection]{Corollary}
\newtheorem{lem}[subsection]{Lemma}
\newtheorem{prop}[subsection]{Proposition}

\newtheorem{thm}[subsection]{Theorem}

\theoremstyle{remark}

\newcommand{\thmref}[1]{Theorem~\ref{#1}}
\newcommand{\secref}[1]{Sect.~\ref{#1}}
\newcommand{\lemref}[1]{Lemma~\ref{#1}}
\newcommand{\propref}[1]{Proposition~\ref{#1}}

\newcommand{\nc}{\newcommand}
\nc{\renc}{\renewcommand}

\nc{\ssec}{\subsection}
\nc{\sssec}{\subsubsection}

\nc{\bi}{\bibitem}
\nc{\on}{\operatorname}

\nc{\A}{\mathcal A}
\nc{\U}{\mathcal U}
\nc{\B}{\mathcal B}
\nc{\C}{\mathcal C}
\renc{\S}{\mathcal S}
\renc{\O}{\mathcal O}
\nc{\M}{\mathcal M}
\nc{\N}{\mathcal N}
\nc{\D}{\mathcal D}
\nc{\K}{\mathcal K}
\nc{\W}{\mathcal W}
\nc{\G}{\mathcal G}
\renc{\L}{\mathcal L}
\renc{\H}{\mathcal H}
\nc{\F}{\mathcal F}

\nc{\wt}{\widetilde}

\nc\tL{\widetilde{L}}

\nc{\CC}{\mathbb C}
\nc{\NN}{\mathbb N}
\nc{\ZZ}{\mathbb Z}
\renc{\AA}{\mathbb A}

\nc{\grg}{\mathfrak g}
\nc\tg{\widetilde{\grg}}
\nc{\grD}{\mathfrak D}
\nc\Hsi{\on{H}^{\frac{\infty}{2}}}
\nc\Hsik{\on{H}^{\frac{\infty}{2}+k}}
\nc\lemb{\mathfrak l}
\nc\remb{\mathfrak r}
\nc\temb{\mathfrak t}
\nc\iemb{\mathfrak i}
\nc{\Hom}{\on{Hom}}

\title{Differential operators on the loop group via chiral algebras} 

\author{S.~Arkhipov and D.~Gaitsgory}

\address{S.A.: Ind. Moscow Univ., 11 Bolshoj Vlasjevskij per., Moscow, 121002 Russia; \newline
 D.G.: Dept. of Math., Harvard University, Cambridge MA 02138}

\email{hippie@mcmme.ru, gaitsgde@math.harvard.edu}

\begin{document}

\maketitle

\section*{$0.$ Introduction}

\ssec{}

Our motivation for writing this paper was twofold. On the one hand, we are trying to make sense
of the notion of D-module on some infinite-dimensional variety attached to an algebraic group,
and on the other hand, we study the structure of a certain interesting representation of the
corresponding affine algebra. 

First, we will explain the representation-theoretic side of the picture.
Thus, let $G$ be an affine algebraic group with Lie algebra $\grg$. 
For simplicity, in the introduction we will assume that $G$ is unimodular. 

Consider the group-scheme $G[[t]]$, which classifies maps of the formal disc 
$\D=\on{Spf}(\CC[[t]])$ into $G$. In addition, given an invariant symmetric form  
$Q:\grg\otimes \grg\to\CC$ one can consider the affine Lie algebra
$\tg_Q:=\grg((t))\oplus \CC$. As a subalgebra it contains $\grg[[t]]$, which is the Lie algebra of $G[[t]]$,
in particular, the ring of regular functions $\O_{G[[t]]}$ is a $\grg[[t]]$-module, via the action
by {\it left-invariant} vector fields.
Let us consider the induced module over the affine algebra ${\mathbb V}_{G,Q}:=\on{Ind}_{\grg[[t]]\oplus
\CC}^{\tg_Q}(\O_{G[[t]]})$.

It turns out that, in addition to being just a $\tg_Q$-module, ${\mathbb V}_{G,Q}$ possesses several other
structures and the goal of this paper is to study them.

\smallskip

The initial observation is that ${\mathbb V}_{G,Q}$ has a natural structure of a {\it vertex operator algebra}.
In this paper we will adopt the language of {\it chiral algebras}, rather than vertex operator algebras,
and the above assertion reads as follows: on every curve $X$ there exists a chiral algebra $\grD_{G,Q}$ and
for a point $x\in X$, there is an isomorphism between ${\mathbb V}_{G,Q}$ and the fiber of $\grD_{G,Q}$ at $x$
for every choice of a local parameter at $x$.

In fact, the $\tg_Q$-module structure on $\grD_{G,Q}$ comes from an embedding 
of chiral algebras $\lemb:\A_{\grg,Q}\to \grD_{G,Q}$, 
where $\A_{\grg,Q}$ is the (Kac-Moody) chiral
algebra canonically attached to the pair $(\grg,Q)$, cf. \secref{KacMoody}.

\smallskip

Here comes a crucial observation, which is the main result of this paper: 

It turns out that in addition to $\lemb$, 
there is another embedding of chiral algebras $\remb:\A_{\grg,Q'}\to \grD_{G,Q}$, where $Q'=-Q-Q_0$, 
where $Q_0$ is the Killing form.

Moreover, the embeddings $\lemb$ and $\remb$ commute with one another in the appropriate sense 
(in the chiral terminology, they *-commute). 
On the level of representation theory, this means that on ${\mathbb V}_{G,Q}$ there is an action of $\tg_{Q'}$, which
commutes with the initial $\tg_Q$-action; the induced action of the subalgebra $\grg[[t]]\subset \tg_{Q'}$ 
comes from the natural $\grg[[t]]$-action on $\O_{G[[t]]}$ by {\it right-invariant} vector fields.

\ssec{}

Now, let us explain the algebro-geometric meaning of the module ${\mathbb V}_{G,Q}$ and of our constructions.

Let $Z$ be an arbitrary affine smooth variety  and let $Z((t))$ be the ind-scheme that classifies maps 
of the formal punctured disc to $Z$. 

One may wonder whether it is possible to define the notion of D-module on $Z((t))$. The general answer in still
not clear, however, a partial solution has been proposed by Beilinson and
Drinfeld:

One can single out a class of chiral algebras, which can be called {\it chiral algebras of differential
operators} (CADO) corresponding to $Z$ (cf. \secref{CDO}). Given a CADO $\grD_Z$, one can consider the category
of chiral modules over it, which will be a candidate for the sought-for category of D-modules. 
For example, the vacuum module, i.e. the fiber of $\grD_Z$ at some point $x\in X$, corresponds to the 
$\delta$-function on the subscheme $Z[[t]]$ inside $Z((t))$. 

We remark that the same class of chiral algebras was independently discovered in \cite{MSV}. 

A priori, there is more than one CADO that corresponds to $Z$ 
(or maybe none) and one thus obtains nonnequivalent
theories of D-modules. Note that even in the finite-dimensional case,
along with ordinary D-modules one can consider twisted D-modules,
and the corresponding categories would not be equivalent in general.
The essential complication of the infinite-dimensional situation 
(which we view in terms of chiral alebras) is that there is 
no preferred (i.e. zero) twisting.

\smallskip

Now let us take $Z=G$. In this case, it turns out that a choice of a 
form $Q:\grg\otimes\grg\to\CC$ 
determines a specific CADO, which will be our $\grD_{G,Q}$. In particular,
in \secref{appendix} we will show that when $Q=0$, the category of chiral 
$\grD_{G,Q}$-modules is indeed close to what one might call the category
of D-modules on $G((t))$.

\medskip

Let us now explain the geometric meaning of the chiral algebra maps
$\lemb$ and $\remb$ mentioned above. The embedding $\lemb$ 
corresponds to the embedding of $\grg((t))$, identified with left-invariant 
vector fields on $G((t))$, into "the ring of differential operators",
and its existence follows fro the construction of $\grD_{G,Q}$.

The embedding $\remb$ corresponds, on the heuristic level, to the 
embedding of $\grg((t))$ into differential operators by means of 
right-invariant vector fields. However, the fact that $\grD_{G,Q}$ is 
attached to $G$ non-canonically results in an anomaly:
in the formula for $\remb$ a non-trivial correction term appears;
in particular, the initial quadratic form $Q$ gets shifted
by the Killing form.
 
\ssec{}

What we said above summarizes the main results and ideas of this paper. Finally, let us mention one more
property of ${\mathbb V}_{G,Q}$, which has to do with the semi-infinite cohomology of $\grg((t))$.

\smallskip

Consider the affine algebra $\tg_{-Q_0}$ corresponding to the Killing form of $\grg$.
Recall that if $N$ is a module over $\tg_{-Q_0}$ (we are assuming that $1\in\CC\subset \tg_{-Q_0}$ acts on $N$
as identity), it makes sense to consider its {\it semi-infinite} cohomology with respect to 
$\grg((t))$, denoted $\Hsik(\grg((t)),N)$.

In particular, if $M$ is a module over $\tg_Q$ for some $Q$, the embedding $\remb$ defines on
the tensor product $N:=M\otimes {\mathbb V}_{G,Q}$ a structure of a module over $\tg_{-Q_0}$, via the diagonal
action. In addition, this tensor product has a commuting $\tg_Q$-module structure by letting
$\tg_Q$ act only on ${\mathbb V}_{G,Q}$ via $\lemb$.

Therefore, the semi-infinite cohomology $\Hsik(\grg((t)),M\otimes {\mathbb V}_{G,Q})$ is naturally a
$\tg_Q$-module. In the last section we prove that, under the assumption that $M$ is $G[[t]]$-integrable,
there is a canonical isomorphism of $\tg_Q$-modules:
$$\Hsik(\grg((t)),M\otimes {\mathbb V}_{G,Q})\simeq M\otimes \on{H}^k(\grg,\CC).$$

Therefore, up to the ordinary Lie algebra cohomology, ${\mathbb V}_{G,Q}$ behaves "as a regular representation"
with respect the operation $M\mapsto \Hsi(\grg((t)),M\otimes {\mathbb V}_{G,Q})$. 

Let us remark in conclusion, that by taking instead
of the vacuum module ${\mathbb V}_{G,Q}$ another chiral module over $\grD_{G,Q}$ (corresponding to the Iwahori
subgroup of $G((t))$, rather than to $G[[t]]$), one could eliminate the appearance of 
$\on{H}^k(\grg,\CC)$ in the above formula. In the terminology suggested
by B.~Feigin, a $\tg_Q$-$\tg_{Q'}$ bimodule with this property should be called {\it a semi-regular
module over the affine algebra}.

\ssec{}

Let us briefly explain how the paper is organized:

In Sect. 1 we review the basics of the theory of chiral algebras. 

In Sect. 2 we introduce the
Beilinson-Drinfeld formalism of chiral differential operators. 

In Sect. 3 we present the construction
of our main object of study, i.e. the chiral algebra $\grD_{G,Q}$. 

In Sect. 4 we prove the main theorem about the existence of the 
embedding $\remb$. 

In Sect. 5 we discuss issues related to semi-infinite cohomology.

In the Appendix (\secref{appendix}) we study the relation between the
category of chiral $\grD_{G,Q}$-modules and ``actual'' D-modules on the
loop group $G((t))$.

\medskip

\ssec{Acknowledgements.}

We would like to express our deep gratitude to B.~Feigin for explaining to us 
his idea that the semi-regular module over the affine algebra should 
possess a structure of vertex algebra. 

All we know about chiral algebras we learned from A.~Beilinson. 
In particular,
our paper uses the unpublished treatese on chiral algebras \cite{BD}, and
in the Appendix we use another unpublished work of Beilinson and
Drinfeld, namely \cite{BD1}. \footnote{Although unpublished,
the above manuscripts are available electronically, cf. our
list of references.}

When this paper was in preparation, we received a preprint by 
V.~Gorbunov, F.~Malikov and V.~Schechtman, 
in which results similar to ours were obtained; we thank 
V.~Schechtman for this communication.

\ssec{Notation and conventions}   \label{notation}

For a scheme $Z$, $\O_X$ (resp., $T_Z$, $\Omega_Z$, $D_Z$) will denote the structure sheaf
(resp., the sheaf of $1$-forms, the sheaf of differential operators) on $Z$.

We will work with a fixed smooth algebraic curve $X$ over $\CC$. For $n\geq 2$, $j_n$ will denote the
(open) embedding of the complement of the diagonal divisor in $X^n$; $\Delta_n$ will denote the
(closed) embedding of the main diagonal $X\to X^n$; the subscript $n$ will sometimes be omitted
when $n=2$ and no confusion is likely to occur.

Our main objects of study will be D-modules on $X$ or $X^n$. Normally, we will work with right
D-modules; for a right D-module $\M$ on $X$, we will denote by $\M^l$ the corresponding left D-module, i.e.
$\M\simeq \M^l\otimes \Omega_X$. For two right D-modules $\M_1$ and $\M_2$, we will denote
by $\M_1\overset{!}\otimes \M_2$ their D-module tensor product, i.e. as an $\O_X$-module
$\M_1\overset{!}\otimes \M_2\simeq \M_1\otimes \M_2\otimes \Omega_X^{-1}$.

\smallskip

\noindent{\bf De Rham cohomology.} Let $\M$ be a right D-module on $X^0$, where $X^0$ is either
an affine curve, or a formal (resp., formal punctured) disc. It then makes sense to consider
the $0$-th De Rham cohomology of $\M$ on $X^0$. By definition, $DR(X^0,\M)$ is a vector space
equal to $\M/\M\cdot T_{X^0}$. We will denote by $h$ the natural projection $h:\Gamma(X^0,\M)\to DR(X^0,\M)$. 

Let $\M$ be a right D-module on $X$ and consider its direct image $\Delta_!(\M)$
under $\Delta:X\to X\times X$. Then the {\it $\O$-module direct image} of $\Delta_!(\M)$ under
$p_2:X\times X\to X$ is naturally a right D-module on $X$ and it surjects by means of the residue map
onto its {\it D-module direct image}, the latter being isomorphic to $\M$ itself. Since this map is
a version of the De Rham projection, we will denote it by $(h\boxtimes \on{id}):\Delta_!(\M)\to \M$.

\section{Chiral algebras}

\ssec{Definition}

We begin by recalling some basic definitions from the theory
of chiral algebras. 

A chiral algebra over $X$ is a right D-module $\A$ endowed with
a {\it chiral bracket}, i.e. a map
$$\{\,\cdot,\cdot\}:j_*j^*(\A\boxtimes\A)\to\Delta_!(\A)$$ which
is antisymmetric and satisfies the Jacobi identity in the following
sense:

\bigskip

Let $\M$ be a D- (or $\O$-) module on $X^n$ equivariant with respect to the action of the symmetric group
$S^n$. (The examples are $\A^{\boxtimes n}$, $j_n{}_*j_n^*(\A^{\boxtimes n})$, $\Delta_n{}_!(\A)$, etc.) For
an element $\sigma\in S_n$, we will denote by the same symbol $\sigma$ its action on the space
$\Gamma(X^n,\M)$ of sections of $\M$.

Let us denote by $\{\{\,\cdot,\cdot\},\cdot\}$ the map
$$j_3{}_*j_3^*(\A\boxtimes\A\boxtimes \A)\overset{\{\,\cdot,\cdot\}\boxtimes\on{id}}\longrightarrow
(\Delta_2\boxtimes \on{id})_!(j_2{}_*j_2^*(\A\boxtimes \A))\overset{\{\,\cdot,\cdot\}}\longrightarrow
(\Delta_2\boxtimes\on{id})_!\circ \Delta_2{}_!(\A)\simeq\Delta_3{}_!(\A).$$

\medskip

We must have:

\begin{itemize}

\item
Let $\sigma$ be the transposition acting on $X\times X$. Then for a (local) section $a$ of $j_*j^*(\A\boxtimes\A)$,
$$\sigma(\{\,\cdot,\cdot\}(a))=-\{\,\cdot,\cdot\} (\sigma(a)).$$

\item
Let $\sigma$ be a cyclic permutation $\sigma$ acting on $X\times X\times X$.  
Then for a (local) section $a$ of $j_3{}_*j_3^*(\A\boxtimes\A\boxtimes \A)$
$$\{\{\,\cdot,\cdot\},\cdot\}(a)+\sigma^{-1}(\{\{\,\cdot,\cdot\},\cdot\}(\sigma(a)))+
\sigma^{-2}(\{\{\,\cdot,\cdot\},\cdot\}(\sigma^2(a)))=0 \in \Delta_!(\A)$$

\end{itemize}

A unit in a chiral algebra $\A$ is a map $\Omega_X\to \A$ such that the composition
$$j_*j^*(\Omega_X\otimes\A)\to j_*j^*(\A\otimes\A)\to \Delta_!(\A)$$ is the canonical map
$j_*j^*(\Omega_X\otimes\A)\to \Delta_!(\Omega_X\overset{!}\otimes \A)\simeq \Delta_!(\A)$.

As in the case of associative algebras, a unit, if it exists, is unique. In all our examples,
chiral algebras will have a unit.

For a point $x\in X$, let $DR(\D_x,\A)$ (resp., $DR(\D^*_x,\A)$) be the De Rham cohomology of $\A$ 
over the formal (resp., formal punctured) disc around $x$. Both $DR(\D_x,\A)$ and $DR(\D^*_x,\A)$
are topological Lie algebras. The latter is sometimes called the local completion of $\A$ at $x$,
denoted $\A_{x,loc}$ and it contains the former as a subalgebra (provided that $\A$ is $X$-flat).

Let $\A_x$ denote the D-module fiber of $\A$ at $x$. Then $\A_x$ is naturally a module over 
$DR(\D^*_x,\A)$, called the {\it vacuum module}. Note that we have a short exact sequence of vector spaces
$$0\to DR(\D_x,\A) \to DR(\D^*_x,\A)\to \A_x\to 0,$$
where the last arrow is given by the action of $DR(\D^*_x,\A)$ on the unit
element $1_x\in \A_x$.

Loosely, speaking, the chiral algebra structure on $\A$ is completely 
determined by the action of $DR(\D^*_x,\A)$ on $\A_x$.

\bigskip

A chiral module over a chiral algebra $\A$ is a right D-module $\M$ endowed with an action map
$\on{act}:j_*j^*(\A\boxtimes \M)\to \Delta_!(\M)$ such that:

\begin{itemize}

\item
Let $\sigma$ be the transposition of the first two factors in $X\times X\times X$. Then for
a (local) section $m$ of $j_3{}_*j_3^*(\A\boxtimes\A\boxtimes\M))$
$$\on{act}\circ (\on{id}\boxtimes \on{act})(a)-\sigma^{-1}(\on{act} \circ(\on{id}\boxtimes \on{act})
(\sigma(a)))=\on{act}\circ (\{\,\cdot,\cdot\}\boxtimes\on{id})(a)\in \Delta_3{}_!(\M).$$

\item
The composition
$j_*j^*(\Omega_X\boxtimes\M)\to j_*j^*(\A\boxtimes\M)\to \Delta_!(\M)$ is the canonical map
$j_*j^*(\Omega_X\boxtimes\M)\to \Delta_!(\M)$.

\end{itemize}

For $x\in X$, let $\M_x$ denote the D-module fiber of $\M$ at $x$. Then $DR(\D^*_x,\A)$ (and hence $DR(\D_x,\A)$)
acts on $\M_x$ in
a natural way and this action is continuous with respect to the (only natural) discrete topology on $\M_x$.

Note, that a chiral module need not be flat as on $\O_X$-module. In fact, we will be particularly
interested in chiral modules supported at a point $x\in X$. If $\M$ is such a module, its D-module fiber
at $x$ is of course $0$, and instead we will denote by $\M_x$ the vector space such that $\M=i_x{}_!(\M_x)$,
where $i_x$ denotes the embedding of the point.

\ssec{Commutative chiral algebras}  \label{com}

Let $\B^l$ be a left D-module endowed with a commutative and associative multiplication map
$\B^l\otimes \B^l\to\B^l$ and a unit $\O_X\to \B^l$, both respecting the 
D-module structure. Such $\B^l$ is a called a (commutative
associative) D-algebra.

Let $\B$ be the corresponding right D-module, i.e. $\B=\B^l\otimes\Omega_X$. We obtain a map 
$$\{\,\cdot,\cdot\}:j_*j^*(\B\boxtimes\B)\to \Delta_!(\B\overset{!}\otimes \B)\to \Delta_!(\B)$$
and we leave it to the reader to check (or, alternatively, consult \cite{BD}) that
it satisfies the axioms of a chiral bracket. 

Moreover, in {\it loc. cit.} it is explained that $\B^l\mapsto \B$ establishes 
an equivalence between the category of (com. assoc.) D-algebras and the subcategory
of chiral algebras for which the chiral bracket $\{\,\cdot,\cdot\}$ factors
as $$j_*j^*(\A\boxtimes\A)\to \Delta_!(\A\overset{!}\otimes \A)\to \Delta_!(\A).$$

Note that the latter condition is equivalent to the fact that the composition
$\A\boxtimes\A\to j_*j^*(\A\boxtimes\A)\overset{\{\,\cdot,\cdot\}}\longrightarrow \Delta_!(\A)$
vanishes. Such chiral algebras will be henceforth referred to as commutative. 
(More generally, for an arbitrary chiral algebra $\A$ one says that
two sections $a_1,a_2\in\A$ {\it *-commute} if $\{\,\cdot,\cdot\}(a_1\boxtimes a_2)=0$.)

\bigskip

The simplest way to produce a D-algebra (and hence a commutative chiral algebra) is as follows:
Let $\M^l$ be a left D-module and let $\B^l$ be the symmetric algebra $\B^l:=\on{Sym}_{\O_X}(\M^l)$. 
We will say that a D-algebra is finitely generated if it is a quotient of a one of 
the above form for a finitely generated D-module $\M^l$.

\medskip

Let $\C$ be a quasi-coherent sheaf of (com. assoc.) algebras on $X$.
We define a D-algebra $J(\C)^l$ as a quotient of 
$\on{Sym}_{\O_X}(D_X\underset{\O_X}\otimes \C)$ by the ideal generated by
$(1\otimes c_1)\cdot (1\otimes c_2)-(1\otimes c_1\cdot c_2)$
and $(1\otimes 1)-1$. This $J(\C)^l$ is called the 
{\it jet construction} of $\C$ and it has the following universal property:
for a D-algebra $\B^l$,
$$\Hom_{D-alg}(J(\C)^l,\B^l)=\Hom_{\O-alg}(\C,\B^l).$$
In other words, $J$ is the left adjoint to the forgetful functor 
from the category of D-algebras
to the category of quasi-coherent sheaves of algebras. 
It is easy to see that $J(\C)^l$ is finitely generated
if $\C$ is.

\medskip

Let $\B^l$ be a D-algebra and $x\in X$ a point (we will denote by $k_x$ the residue field of $X$ at $x$).
Let $\B_x$ be the D-module fiber of $\B^l$ at $x$.
By definition, as a set $\on{Spec}(\B_x)$ consists of all $\O_x$-algebra homomorphisms $\B^l\to k_x$.
It is well-known that for a left D-module $\M^l$ on $X$, 
$$\Hom_{\O_X}(\M^l,k_x)=\Hom_{D_X}(\M^l,\hat{\O}_x).$$
Analogously, we have:
$$\Hom_{\O_X}(\B^l,k_x)=\Hom_{D_X}(\B^l,\hat{\O}_x),$$
as algebras. Therefore, $\on{Spec}(\B_x)$ consists of all flat sections of $\on{Spec}(\B^l)$ over $\D_x$. 

In particular, if $\B^l=J(\C)^l$, by the universal property, 
$\on{Spec}(J(\C)_x)$ consists of all sections of $\C$ over $\D_x$, which explains the name ``jets''.

\ssec{Lie-* algebras}  \label{Liestar}

To proceed further, we need to recall the definition of Lie-* algebras. These objects are basic
tools for constructing non-commutative chiral algebras.

By definition, a Lie-* algebra is a right D-module $L$ on $X$ endowed with a Lie-* bracket
$\{\,\cdot,\cdot\}:L\boxtimes L \to\Delta_!(L)$, which is antisymmetric and satisfies the Jacobi
identity in the same sense as in the definition of a chiral bracket.

A Lie-* algebra is called commutative if its Lie-* bracket is trivial, i.e. equals $0$.

\medskip

A Lie-* module over a Lie-* algebra $L$ is a right D-module $\M$ endowed with a map:
$$\on{act}:L\boxtimes \M\to\Delta_!(\M),$$
which satisfies the condition identical to condition (1) in the definition of chiral modules
over a chiral algebra.

If $L$ is a Lie-* algebra, the vector spaces $DR(\D_x,L)$ and $DR(\D^*_x,L)$ are naturally topological
Lie algebras. If $\M$ is a Lie-* module over $L$, then its fiber $\M_x$ is a continuous module over
$DR(\D_x,L)$. 

\medskip

Here is a simplest example of a Lie-* algebra. Let $\grg$ be a Lie algebra and let us consider 
$L_{\grg}=\grg\otimes D_X$.  Then the Lie bracket on $\grg$ induces a map
$$\grg\otimes D_X\boxtimes \grg\otimes D_X\to \Delta_!(\grg\otimes D_X),$$
which is easily seen to satisfy the properties of a Lie-* bracket.

\ssec{The universal enveloping chiral algebra}   \label{envelop}

There exists an obvious forgetful functor from the category of chiral algebras to that of Lie-* algebras:
we compose the chiral bracket $\{\,\cdot,\cdot\}:j_*j^*(\A\boxtimes\A) \to\Delta_!(\A)$ with the embedding
$\A\boxtimes \A\to j_*j^*(\A\boxtimes\A)$. (Note that the resulting Lie-* algebra is commutative if and only if
$\A$ is a commutative chiral algebra.)

It is a basic fact that the above functor admits a left adjoint $L\mapsto \U(L)$. In other words, for
a chiral algebra $\A$ there is a functorial isomorphism
$$\Hom_{Lie-*}(L,\A)=\Hom_{Chiral}(\U(L),\A).$$
This $\U(L)$ is called the {\it universal enveloping} chiral algebra of $L$.

We refer the reader to \cite{BD} or \cite{Ga} for the proof of the existence of this functor and
of its basic properties, some of which are reviewed below.
In what follows we will assume that $L$ is flat as an $\O_X$-module.

\medskip

By the universal property, we have a map $L\to \U(L)$. Hence, if $\M$ is a chiral module
over $\U(L)$ supported at $x\in X$, then $\M_x$ is a continuous module not only over 
$DR(\D_x,L)$, but over the whole $DR(\D^*_x,L)$.

\begin{lem}  \label{catmodules}
The above functor establishes an equivalence
$$\text{Chiral modules over $\U(L)$ supported at $x$ } \leftrightarrow \text{ Continuous modules over $DR(\D^*_x,L)$.}$$
\end{lem}

In particular, $\U(L)_x$ is a $DR(\D^*_x,L)$-module.

\begin{lem} \label{fiber}
We have:
$$\U(L)_x\simeq \on{Ind}^{DR(\D^*_x,L)}_{DR(\D_x,L)}(\CC),$$
where $\CC$ is the trivial module over $DR(\D_x,L)$, in such a way that

\smallskip

\noindent {\em a)} The generator in $\on{Ind}^{DR(\D^*_x,L)}_{DR(\D_x,L)}(\CC)$ corresponds to the unit 
$1_x\in \U(L)_x$. 

\smallskip

\noindent {\em b)}
The embedding $L\to \U(L)$ corresponds 
to $$L_x=DR(\D^*_x,L)/DR(\D_x,L)\to \on{Ind}^{DR(\D^*_x,L)}_{DR(\D_x,L)}(\CC).$$
\end{lem}

\medskip

Here is an additional property of the chiral algebra $\U(L)$ (\cite{BD,Ga}). 

\begin{prop}  \label{filtration}
As a D-module, $\U(L)$ carries a unique increasing filtration
$\U(L)=\underset{i\geq 0}\cup \U(L)_i$ such that 

\smallskip 

\noindent {\em a)}
$\U(L)_0=\Omega_X$

\smallskip 

\noindent {\em b)}
The embedding $L\to \U(L)$ induces an isomorphism 
$\U(L)_1\simeq\Omega_X\oplus L$ 

\smallskip 

\noindent {\em c)}
It is compatible with the chiral bracket in the sense that
\begin{align*}
&\{\,\cdot,\cdot\}:j_*j^*(\U(L)_i\boxtimes \U(L)_j)\mapsto \Delta_!(\U(L)_{i+j})  \\
&\{\,\cdot,\cdot\}:\U(L)_i\boxtimes \U(L)_j\mapsto \Delta_!(\U(L)_{i+j-1}).
\end{align*}

\end{prop}

In particular, the above proposition 
implies that $\on{gr}(\U(L)):=\underset{i\geq 0}\oplus \U(L)_j/\U(L)_{j-1}$ is a 
commutative chiral algebra. The map $L\to \U(L)$ induces a map $\on{Sym}(L)\to \on{gr}(\U(L))$
and it is shown in \cite{BD} that an analog of the PBW theorem holds:

\begin{thm}  \label{chiralPBW}
$\on{Sym}(L)\simeq \on{gr}(\U(L))$.
\end{thm}

On the level of fibers, this filtration corresponds to the standard 
filtration on the induced module 
$$U(DR(\D^*_x,L))\underset{U(DR(\D_x,L))}
\otimes\CC\simeq\on{Ind}^{DR(\D^*_x,L)}_{DR(\D_x,L)}(\CC)\simeq \U(L)_x,$$
coming from the filtration on the usual universal enveloping 
algebra $U(DR(\D^*_x,L))$.

\section{Chiral differential operators}  \label{CDO}

\ssec{The case of an affine space}   \label{flatcase}

As was explained in the introduction, the object of study of this paper is a 
chiral algebra $\grD_{G,Q}$ attached to a group $G$ and a form $Q$ on $\grg$. This chiral algebra
will be constructed as a {\it chiral algebra of differential operators} on $G$. In this section we 
will recall the theory of chiral differential operators, developed by 
Beilinson-Drinfeld and Malikov-Schechtman-Vaintrob, cf. \cite{BD},
Sect. 3.9 and \cite{MSV}.

Let $Z$ be a smooth affine algebraic variety and let $\C$ be the 
$\O_X$-algebra 
$\O_Z\otimes\O_X$. Consider $Z(\hat\O_x):=\on{Spec}(J(\C)^l_x)$. 
This is a scheme, which is
a projective limit of finite-dimensional schemes $Z_i:=Z(\hat\O_x/m_x^i)$, 
where $m_x\subset \hat\O_x$
is the maximal ideal.

Since the reduction maps $Z(\hat\O_x/m_x^i)\to Z(\hat\O_x/m_x^{i-1})$ are smooth, one can easily make
sense of the category of {\it left} D-modules on $Z(\hat\O_x)$: by definition, this category 
is $\underset{\longrightarrow}{lim} \,D$-{mod}$(Z_i)$. 

In addition, one can define the ind-scheme $Z(\hat\K_x)$, 
where $\K_x$ is the fraction field of $\O_x$
(cf. \cite{Ga1}). In general, if 
one has an ind-scheme of ind-finite type (i.e. representable
as a union of closed subschemes, each of which is of finite type) 
one can define
the notion of a {\it right} D-module on it. However, $Z(\hat\K_x)$ 
is not of ind-finite type and it is 
a priori not clear that the theory of D-modules on it exists at all. 

Thus, the first obstacle in defining D-modules on $Z(\hat\K_x)$ is 
the fact that one cannot pass between left and right D-modules in 
the infinite-dimensional setting. Therefore, we will start by
considering the example of $Z=\AA^n$, in which case this difficulty can 
be eliminated, and which will suggest the form of the solution in general. 

We emphasize that material that we are reviewing in this section 
is fully contained in \cite{BD} and in \cite{MSV}.

\bigskip

Recall that if $V$ is a finite-dimensional vector space, then the category 
of left (resp., right) D-modules on $V$
viewed as an algebraic variety is naturally equivalent to the category of 
left (resp., right) modules
over the Weyl algebra, $W(V)$. By definition, $W(V)$ is generated by 
elements $v\in V$, $v^*\in V^*$ which
satisfy the relations
\begin{align*}
& v_1\cdot v_2=v_2\cdot v_1,\,\,\, v^*_1\cdot v^*_2=v^*_2\cdot v^*_1  \\
& v\cdot v^*-v^*\cdot v=\langle v,v^*\rangle.
\end{align*}

This definition can be generalized to the case when $V$ is no longer finite-dimensional, but rather is a topological
vector space, e.g. $V:=\AA^n\otimes \hat\K_x$. Note that in the latter case, the topological dual of $V$
can be identified with the space of $1$-forms on $\D^*_x$ with values in $\AA^n$, i.e. with
$\AA^n\otimes \Omega_{\hat\K_x}$.

Now, we claim that the category of
continuous left modules over $W(\AA^n\otimes \hat\K_x)$ 
is naturally equivalent to
to the category of modules supported at $x$ over an explicit 
chiral algebra, called the Weyl algebra.

\bigskip

Consider the Heisenberg Lie-* algebra
$$\H^n=\Omega_X\oplus D_X^{\oplus n} \oplus (\Omega_X\otimes D_X)^{\oplus n},$$
where the only non-trivial component of the Lie-* bracket,
\begin{align*}
&D_X^{\oplus n}\boxtimes (\Omega_X\otimes D_X)^{\oplus n}\to \Delta_!(\Omega_X) \text{ and }\\
&(\Omega_X\otimes D_X)^{\oplus n}\boxtimes D_X^{\oplus n}\to \Delta_!(\Omega_X),
\end{align*}
are given by means of the canonical element in
$$\Hom_{D_X}(D_X\boxtimes (\Omega_X\otimes D_X),\Delta_!(\Omega_X))\simeq 
\Hom_{\O_X}(\O_X\boxtimes \Omega_X,\Delta_!(\Omega_X)).$$

Note that $DR(\D^*_x,\H^n)$ is the usual Heisenberg Lie algebra 
$$(\AA^n\otimes\hat\K_x)\oplus (\AA^n\otimes\Omega_{\hat\K_x})\oplus \CC,$$ where
$DR(\D^*_x,\Omega_X)$ is identified with $\CC$ via the residue map.

Consider the universal enveloping chiral algebra $\U(\H^n)$. We have two different embeddings of
the D-module $\Omega_X$ into $\U(\H^n)$. One is the unit in $\U(\H^n)$ and another comes from $\Omega_X\subset \H^n$.
Let $I$ be the ideal in $\U(\H^n)$ generated by the anti-diagonal copy of $\Omega_X$.

Set $\W^n=\U(\H^n)/I$. We have the following:

\begin{prop}
The category of chiral $\W^n$-modules supported at $x$, is equivalent to the category
of modules over the Weyl algebra $W(\AA^n\otimes\hat\K_x)$.
\end{prop}

\begin{proof}
Using \lemref{catmodules}, we obtain that the category
of $\W^n$-modules supported at $x$, is equivalent to the category of
continuous modules over $DR(\D^*_x,\H^n)$, on which
$1\in \CC$ acts as the identity. However, these are the same as modules over $W(\AA^n\otimes\hat\K_x)$.
\end{proof}

To conclude the discussion of the ``flat'' case (i.e. $Z=\AA^n$), let us make the following observation:
$\Omega_X\otimes D_X^{\oplus n}$ is contained in $\H^n$ as a Lie-* subalgebra; hence, we obtain
a map of chiral algebras $\U(\Omega_X\otimes D_X^{\oplus n})\to \W^n$. However, since $\Omega_X\otimes D_X^{\oplus n}$
is commutative, $\U(\Omega_X\otimes D_X^{\oplus n})\simeq \on{Sym}(D_X^{\oplus n})\otimes\Omega_X$, 
the latter being isomorphic to the jet construction of 
$\C=\O_{\AA^n}\otimes\O_X$. In addition, $\H^n$ contains as a Lie-* subalgebra $D_X^{\oplus n}$,
and $\W^n$ is generated by $D_X^{\oplus n}$ and 
$J(\on{Sym}(\O_X^{\oplus n}))$. This will be the form of the answer in
general.

\ssec{The Lie-* algebra of vector fields}  \label{vectorfields}

Let $\C$ be a locally finitely generated quasi-coherent sheaf of algebras over $X$,
such that $\on{Spec}(\C)$ is smooth over $X$. (In practice, we will take $\C=\O_Z\otimes\O_X$.)
Consider $\B^l=J(\C)^l$. In principle, our discussion applies to a more general D-algebra,
but for simplicity, we will consider the above case only.

Let $T_\C$ be the $\C$-module of vertical vector fields on $\on{Spec}(\C)$. 
Consider the tensor product
$$\Theta_\C:=T_\C\underset{\C}\otimes  (J(\C)^l\underset{\O_X}\otimes D_X).$$

This is a right D-module on $X$ and it carries an action of the D-algebra $J(\C)^l$. In particular,
we obtain a map
$$j_*j^*(J(\C)\boxtimes \Theta_\C)\to \Delta_!(J(\C)\overset{!}\otimes \Theta_\C)\to \Delta_!(\Theta_\C),$$
which makes it a module over $J(\C)$, viewed as a chiral algebra.

In addition, we claim that $\Theta_\C$ has a natural structure of a Lie-* algebra over which
$J(\C)$ is a module:

\begin{lem}   \label{vector kitchen}
There exists a unique Lie-* bracket $\{\,\cdot,\cdot\}:\Theta_\C\boxtimes \Theta_\C\to\Delta_!(\Theta_\C)$
and an action map $\on{act}:\Theta_\C\boxtimes J(\C)\to \Delta_!(J(\C))$ such that

\smallskip

\noindent {\em a)} The induced maps
$$T_\C\boxtimes T_\C\to \Theta_\C\boxtimes \Theta_\C\to\Delta_!(\Theta_\C)
\text{ and } T_\C\boxtimes (\C\otimes\Omega_X)\to \Theta_\C\boxtimes J(\C) \to\Delta_!(J(\C))$$
factor, respectively, through the natural maps
$$T_\C\boxtimes T_\C\overset{\text{Lie bracket}}\longrightarrow T_\C\hookrightarrow \Delta_!(\Theta_\C)\text{ and }
T_\C\boxtimes (\C\otimes\Omega_X)\overset{\text{Lie der.}}
\longrightarrow \C\otimes\Omega_X\hookrightarrow \Delta_!(J(\C)).$$

\smallskip

\noindent {\em b)}
They satisfy the Leibnitz rule in the natural sense.
\end{lem}

Let us try to give an intuitive view point on $\Theta_\C$. Consider a point in the fiber of
$\on{Spec}(J(\C)^l)$ over $x\in X$. As was explained in \secref{com}, it corresponds to a
section $\phi':\D_x\to \on{Spec}(\C)$, which is the same as
a flat section $\phi:\D_x\to \on{Spec}(J(\C)^l)$. We can view
$\Theta_\C$ as a coherent sheaf with a right $D_X$-action on $\on{Spec}(J(\C)^l)$ 
and let us consider its pull-back with respect to
$\phi$ as a right D-module on $\D_x$. 

\begin{lem}
There is a natural isomorphism $DR(\D_x,\phi^*(\Theta_\C))\simeq \Gamma(\D_x,\phi'{}^*(T_\C))$.
\end{lem}

Note that the RHS of the isomorphism of the lemma is canonically isomorphic to the vertical tangent space
to $\on{Spec}(J(\C)^l)$ at our chosen point in the spectrum.

\ssec{The Beilinson-Drinfeld definition of chiral differential operators}

For a smooth algebraic variety $Z$, let $D_Z$ denote the algebra of differential operators.
By definition, $D_Z$ carries a canonical filtration $D_Z=\underset{i\geq 0}\cup D_Z^i$ such that

\begin{enumerate}

\item 
For $a\in D_Z^i$, $b\in \D_Z^j$ we have:
$$a\cdot b\in D_Z^{i+j},\,\,\, [a,b]\in D_Z^{i+j-1}.$$

\item
$\D_Z^0=\O_Z$ as algebras

\item
We have an isomorphism $D_Z^1/D_Z^0\simeq T_Z$, which respects a) the Lie algebra structure,
b) the structure of an $\O_Z$-module on $T_Z$, c) The adjoint action of $T_Z$ on $\O_Z$.

\item
We have a splitting $T_Z\to D_Z^1$, which respects the Lie-algebra and the $\O_Z$-module
structure, where on $D_Z^1$, $\O_Z$ acts by {\it left} multiplication.

\item
The natural map $\on{Sym}_{\O_Z}(T_Z)\to \on{gr}(D_Z):=\oplus D_Z^i/D_Z^{i-1}$ is an isomorphism.

\end{enumerate}

\medskip

In \cite{BB} Beilinson and Bernstein introduced the notion of twisted differential operators (TDO).
By definition, a TDO on $Z$ is a quasi-coherent sheaf $D'_Z$ of associative algebras over $Z$
endowed with an increasing filtration $D'_Z=\underset{i\geq 0}\cup (D'_Z)^i$ such that
properties (1), (2), (3) and (5) hold. In particular, the untwisted differential operators
are distinguished by the embedding of point (4) above.

\medskip

The idea of the approach of Beilinson and Drinfeld is
construct a suitable chiral algebra $\grD_Z$, 
which we will call the chiral algebra of differential
operators (CADO),
and to {\it define} D-modules on $Z(\hat\K_x)$
as chiral $\grD_Z$-modules supported at $x\in X$.

When $Z\simeq\AA^n$, the Weyl chiral
algebra $\W^n$ will be an example of a CADO. In the Appendix, \secref{finale}
we will show that when $Z$ is an algebraic group, the category of
chiral $\grD_Z$-modules is indeed a reasonable approximation to the
ill-defined category of D-modules on $G(\hat\K_x)$.

However, it turns out that there is no CADO canonically associated 
with a variety $Z$. In particular,
it does not always exist, and when it does, it is not unique. 
Moreover, the category of possible CADOs 
on $Z$ forms a gerbe over an explicit Picard category. 

The need for $\grD_Z$ was independently realized
by the authors of \cite{MSV}, where this object was studied in the 
language of vertex algebras, under the name ``chiral structure sheaf''. 
In particular, the authors of \cite{MSV} calculated the obstruction 
for the existence of a CADO on a given $Z$.

\medskip

Here is the definition. Let $\C$ be as in the previous subsection.
A chiral algebra $\grD_\C$ is called a CADO on $\on{Spec}(\C)$ if 
as a D-module on $X$ it is endowed with an increasing filtration 
$\grD_\C=\underset{i\geq 0}\cup \grD_\C^i$ with the following properties:

\begin{enumerate}

\item
\begin{align*}
&\{\,\cdot,\cdot\}:j_*j^*(\grD_\C^i\boxtimes \grD_\C^j)\mapsto \Delta_!(\grD_\C^{i+j})  \\
&\{\,\cdot,\cdot\}:\grD_\C^i\boxtimes \grD_\C^j\mapsto \Delta_!(\grD_\C^{i+j-1}).
\end{align*}

\item
$\grD_\C^0\simeq J(\C)$ as chiral algebras.

\item
We have an homomorphism of $\Theta_\C\to \grD^1_\C/\grD^0_\C$, which respects
a) the Lie-* algebra structure on both sides, b) the structure of a chiral $J(\C)$-module
on both sides, c) The structure on $J(\C)$ of a Lie-* module over $\Theta_\C$ and $\grD_\C^1/\grD^0_\C$.

\item
The natural map of D-algebras 
$\on{Sym}_{J(\C)^l}(\Theta_\C^l)\to \on{gr}(\grD_\C)^l:=\oplus (\grD_\C^i)^l/(\grD_\C^{i-1})^l$
is an isomorphism. (The superscript ``$l$'' denotes, as usual, the corresponding left D-module on $X$.)

\end{enumerate}

The reason for the non-existence of a canonical CADO is that one cannot require the existence
of an embedding $\Theta_\C\to\grD_\C^1$ which is compatible with the structure of a chiral
$J(\C)$-module. The existence of such an embedding would contradict other axioms.

\medskip

To conclude this section, let us go back to the example of $Z=\AA^n$ and $\C=\O_Z\otimes \O_X$.
The construction of the canonical CADO $\W^n$ in this case 
is explained by the fact that we demanded an extra rigidity. We have:
$$\Theta_\C\simeq J(\C)^l\otimes  D_X^{\oplus n}$$
and it contains $D_X^{\oplus n}$ as a Lie-* subalgebra. Our extra rigidity is the embedding of
$D_X^{\oplus n}$ into $\grD^1_\C:=\W^n$ as a Lie-* algebra.

In our case of interest, when $\AA^n$ is replaced
by an arbitrary algebraic group, we will use a similar idea. We should warn the reader that the filtration
on $\W^n$ which comes from the universal enveloping algebra construction is
{\it different} from the canonical CADO filtration. 

\section{The case of an algebraic group}

\ssec{Affine Lie-* algebras}         \label{KacMoody}

Let $G$ be a linear algebraic group and let $\grg$ denote its Lie algebra. In \secref{Liestar} we showed that the
D-module $L_{\grg}:=\grg\otimes D_X$ has a natural structure of a Lie-* algebra.

Let now $Q:\grg\otimes\grg\to\CC$ be a $G$-invariant symmetric form. Consider the D-module
$\tL_{\grg,Q}=\Omega_X\oplus L_{\grg}$. We define on it a Lie-* bracket with the following components:
\begin{align*}
& L_{\grg}\boxtimes L_{\grg}\to \Delta_!(L_{\grg}) \text{ is as before} \\
& \Omega_X\boxtimes L_{\grg,Q}\to \Delta_!(L_{\grg,Q}) \text{ vanishes, and}  \\
& L_{\grg}\boxtimes L_{\grg}\to \Delta_!(\Omega_X) \text{ is defined as follows:}
\end{align*}

To specify such a map is the same as to give a map $\grg\otimes\grg\to\Gamma(X\times X,\Delta_!(\Omega_X))$.
Let $1'\in \Gamma(X\times X,\Delta_!(\Omega_X))$ be the canonical section, killed by the square of the equation
of the diagonal. (In coordinates, $1'$ is the image of $\frac{dx\boxtimes dy}{(x-y)^2}$ under
$j_*j^*(\Omega_X\boxtimes\Omega_X)\to \Delta_!(\Omega_X)$). We set
$$\grg\otimes\grg\overset{Q}\to\CC\overset{1\to 1'}\longrightarrow\Gamma(X\times X,\Delta_!(\Omega_X)).$$

Consider the universal enveloping chiral algebra $\U(\tL_{\grg,Q})$. We define the affine chiral algebra
$\A_{\grg,Q}$ as the quotient of $\U(\tL_{\grg,Q})$ by the ideal generated by the antidiagonal copy of 
$\Omega_X\to \U(\tL_{\grg,Q})$, as in the definition of $\W^n$, cf. \secref{flatcase}.

\medskip

Note that $\tg_Q:=DR(\D_x^*,\tL_{\grg,Q})\simeq (\grg\otimes\hat\K_x)\oplus\CC$ is the standard affine Lie
algebra $\grg\otimes\hat\K_x$ corresponding to $Q$. In what follows, by a representation of $\tg_Q$ we will
mean a continuous representation, on which $1\in\CC\subset \tg_Q$ acts as identity. From \lemref{catmodules} we obtain

\begin{lem} \label{cataffine}
The category of chiral $\A_{\grg,Q}$-modules concentrated at $x$ is naturally equivalent to the category 
of continuous representations of $\tg_Q$. 
\end{lem}

\ssec{Construction of the CADO}   \label{secconstr}
From now on we will specialize to the case $Z=G$. 
We take $\C=\O_G\otimes \O_X$;
however, we shall slightly abuse the notation and will write 
$J(G)$ (resp., $\Theta_G$, $\grD_G$) instead of $J(\C)$ (resp., $\Theta_\C$, $\grD_\C$).

First, let us observe that we have the map $\grg\to T_G$ which corresponds to left invariant vector fields
on $G$. It is easy to see that it extends to a map of Lie-* algebras $L_{\grg}\to \Theta_G$. Moreover,
the above it induces an isomorphism
$$J(G)^l\otimes L_{\grg}\to \Theta_G.$$

\begin{thm} \label{constr}

To a fixed form $Q$ as above, there corresponds a canonical CADO, $\grD_{G,Q}$, with the following extra structure:
there exists an embedding of Lie-* algebras
$$\lemb:\tL_{\grg,Q}\to \grD^1_{G,Q},$$
such that $\Omega_X\subset\tL_{\grg,Q}$ maps identically to $\Omega_X\subset J(G)=\grD^0_{G,Q}$
and the composition 
$$\tL_{\grg,Q}\to \grD^1_{G,Q}/\grD^0_{G,Q}\simeq \Theta_G$$
equals the above canonical map $\tL_{\grg,Q}\to L_{\grg}\to \Theta_G$.

\end{thm}

Let us first explain the idea of the proof. The D-module $L_{\grg}\subset \grD_{G,Q}$ should be thought
of as corresponding to ``left-invariant vector fields on $G(\hat\K_x)$''. Hence our task is to reconstruct
the ``ring of differential operators'' by knowing functions and left-invariant vector fields. 

In the finite-dimensional situation we would proceed as follows: the direct sum $\O_G\oplus \grg$ has a natural
structure of a Lie algebra.
Consider the universal enveloping algebra $U(\O_G\oplus \grg)$. As a subalgebra
it contains $U(\O_G)\simeq \on{Sym}(\O_G)$. However, since $\O_G$ is already an algebra, we have a natural map
$\on{Sym}(\O_G)\to\O_G$. Let $I$ denote the ideal generated by the kernel of this map inside $U(\O_G\oplus \grg)$.

We have: $D_G=U(\O_G\oplus \grg)/I$. In the chiral setting the idea of the proof is exactly the same.

\begin{proof}

Consider the D-module $J(G)\oplus L_{\grg}$. We claim that it has a natural Lie-* algebra structure, whose 
components are as follows:

\begin{align*}
& L_{\grg}\boxtimes L_{\grg}\to \Delta_!(L_{\grg}) \text{ is the old bracket on $L_{\grg}$}  \\
& J(G)\boxtimes J(G)\to \Delta_!(J(G)) \text{ vanishes} \\
& L_{\grg}\boxtimes L_{\grg}\to \Delta_!(J(G)) \text{ equals } 
L_{\grg}\boxtimes L_{\grg}\overset{Q}\to\Delta_!(\Omega_X)\to
\Delta_!(J(G)) \\
& L_{\grg}\boxtimes J(G)\to \Delta_!(J(G)) \text{ is the action map introduced in \secref{vectorfields}.}
\end{align*}

Consider the universal enveloping chiral algebra $\U(J(G)\oplus L_{\grg})$. It contains as a chiral subalgebra
$\U(J(G))$, where $J(G)$ is regarded as a Lie-* algebra with a trivial bracket. However, since $J(G)$
is already a chiral algebra, we have a homomorphism 
$$\U(J(G))\to J(G).$$

Let $I$ denote the ideal inside $\U(J(G)\oplus L_{\grg})$ generated by the kernel of this map. We set
$$\grD_{G,Q}:=\U(J(G)\oplus L_{\grg})/I$$
and we claim that it satisfies all the requirements.

We define the filtration on $\grD_{G,Q}$ by declaring that $\grD^i_{G,Q}$ is the image under the chiral bracket
of $j_*j^*(J(G)\boxtimes \U(J(G)\oplus L_{\grg})_i)$, where $\U(J(G)\oplus L_{\grg})_i$ is the corresponding
term of the filtration of the universal enveloping chiral algebra, cf. \secref{envelop}.

From this definition it is easy to see that the embedding 
$J(G)\to \grD_{G,Q}$ induces
an isomorphism $J(G)\simeq\grD^0_{G,Q}$. Property 1 of the definition of CADO holds in view
of \propref{filtration}(c), since $J(G)$ is commutative. 

The chiral bracket induces a map
$$j_*j^*(L_{\grg}\boxtimes J(G))\to \Delta_!(\grD^1_{G,Q}).$$
When we compose it with the projection $\grD^1_{G,Q}\to \grD^1_{G,Q}/\grD^0_{G,Q}$, we obtain a map,
which factors as
$$j_*j^*(L_{\grg}\boxtimes J(G))\to \Delta_!(L_{\grg}\overset{!}\otimes J(G))\to 
\Delta_!(\grD^1_{G,Q}/\grD^0_{G,Q}).$$

However, $L_{\grg}\overset{!}\otimes J(G)\simeq \Theta_G$, hence we obtain a map
$\Theta_G\to \grD^1_{G,Q}/\grD^0_{G,Q}$, which is easily seen to satisfy the three conditions of point 3. 

\medskip

The map $$\on{Sym}_{J(G)^l}(\Theta_G^l)\to \on{gr}(\grD_\C)^l$$
is a surjection, as follows from \thmref{chiralPBW}. Since the LHS is $X$-flat, to prove that this
map is an isomorphism, it is enough to do so on the level of fibers.

Let $\O_{G(\hat\O_x)}$
denote the ring of functions on the group-scheme $G(\hat\O_x)$. It follows from \lemref{fiber} that
$$(\grD_{G,Q})_x=\on{Ind}_{\grg\otimes \hat\O_x}^{\tg_Q}(\O_{G(\hat\O_x)}),$$
where $\O_{G(\hat\O_x)}$ is a $\grg\otimes \hat\O_x$-module via the action by left-invariant vector fields.
(In the above formula the induction is 
$U'(\tg_Q)\underset{U(\grg\otimes \hat\O_x)}\otimes \O_{G(\hat\O_x)}$, where
$U'(\tg_Q)$ is the quotient of the universal enveloping algebra by the standard
relation that $1\in\CC\subset\tg_Q$ equals the identity.)

\smallskip

The filtration induced on $(\grD_{G,Q})_x$ coincides with the standard filtration on the induced module.
Therefore, as in \thmref{chiralPBW},
$$\on{gr}(\grD_{G,Q})_x\simeq 
\on{gr}(\on{Ind}_{\grg\otimes \hat\O_x}^{\tg_Q}(\O_{G(\hat\O_x)}))\simeq
\O_{G(\hat\O_x)}\otimes \on{Sym}(\grg\otimes\hat\K_x/\grg\otimes \hat\O_x).$$
Therefore, the map
$$\on{Sym}_{J(G)_x}(\Theta_G^l)_x\simeq
\on{Sym}_{\O_{G(\hat\O_x)}}((\grg\otimes\hat\K_x/\grg\otimes \hat\O_x)\otimes \O_{G(\hat\O_x)})\to 
\on{gr}(\grD_{G,Q})_x$$
is an isomorphism.

\end{proof}

In the course of the proof we have shown the following

\begin{cor} \label{semidual}
We have an isomorphism of $\widetilde{\grg}_Q$-modules:
$${\mathbb V}_{G,Q}:=(\grD_{G,Q})_x\simeq\on{Ind}_{\grg\otimes \hat\O_x}^{\tg_Q}(\O_{G(\hat\O_x)}).$$
\end{cor}

\ssec{The main result}   \label{mainresult}

As we have explained before, the D-module $L_{\grg}$ embedded into $\grD_{G,Q}$ corresponds
to left-invariant vector fields on the loop group. It is, therefore, natural to ask whether there exists
another embedding $\remb:L_{\grg}\to \grD_{G,Q}$, which *-commutes with $\lemb(L_{\grg})$ and corresponds
to right-invariant vector fields. In particular, one would expect the module
$\on{Ind}_{\grg\otimes \hat\O_x}^{\tg_Q}(\O_{G(\hat\O_x)})$, which corresponds
to the ``$\delta$-function on $G(\hat\O_x)$ inside $G(\hat\K_x)$'' to be in fact a bi-module over
the affine algebra.

The answer to this question is positive. However, this second embedding of $L_{\grg}$ develops an
anomaly. The precise statement is given below.

\bigskip

First, let us introduce the following canonical extension $\F$ of $\O_X$-modules on $X$:
$$0\to \Omega_X\to\F\to\O_X\to 0.$$
We take $\F:={p_2}_*(\O_X\boxtimes \Omega_X(\Delta)/\O_X\boxtimes \Omega_X(-\Delta))$, where $p_2:X\times X\to X$ is the
second projection. Note that this extension is globally non-trivial, i.e. when $X$ is compact its class is the canonical
element in $H^1(X,\Omega_X)$.

Let $\grg$ be a Lie algebra and let $\rho:\grg\to\CC$ be its modular character (i.e. the character by which $\grg$
acts on $\Lambda^{\on{top}}(\grg^*)$). We define the D-module 
$(\grg\otimes D_X)'$ as the quotient of
$\grg\otimes (\F\otimes D_X)$ by the kernel of the map 
$$\grg\otimes (\Omega_X\otimes D_X)\overset{\rho}\to \Omega_X\otimes D_X\to D_X.$$
By construction, $(\grg\otimes D_X)'$ is a D-module extension
$$0\to \Omega_X\to (\grg\otimes D_X)'\to \grg\otimes D_X\to 0$$ and we have a splitting
$[\grg,\grg]\otimes D_X\to (\grg\otimes D_X)'$, since $\rho$ is trivial on $[\grg,\grg]$.

Given an invariant form $Q$ on $\grg$ we obtain a Lie-* algebra structure on $(\grg\otimes D_X)'$, 
which we will denote by $\tL'_{\grg,Q}$, defined in the same way as in the case of $\tL_{\grg,Q}$.
Let $Q_0:\grg\otimes\grg\to\CC$ be the Killing form. We define the involution 
$Q\to Q'$ on the set of invariant quadratic forms by $Q'=-Q_0-Q$.

\begin{thm} \label{main}
Let $\grD_{G,Q}$ be as in \thmref{constr}. There exists a canonical Lie-* algebra embedding
$\remb:\tL'_{\grg,Q'}\to \grD_{G,Q}$ with the following properties:

\smallskip

\noindent {\em a)}
$\Omega_X\subset \tL'_{\grg,Q'}$ maps identically onto $\Omega_X\subset \grD_{G,Q}$.

\smallskip

\noindent {\em b)}
The composition
$$\tL_{\grg,Q}\boxtimes \tL'_{\grg,Q'}\overset{\lemb\boxtimes \remb}\longrightarrow \grD_{G,Q}\boxtimes \grD_{G,Q}
\overset{\{\,\cdot,\cdot\}}\longrightarrow \Delta_!(\grD_{G,Q})$$
vanishes. 

\smallskip

\noindent {\em c)}
The image of $\remb$ lies in $\grD^1_{G,Q}$ and the resulting map 
$$L_{\grg}\to \grD^1_{G,Q}/\grD^0_{G,Q}\simeq \Theta_G$$
corresponds to the embedding of $\grg$ into $T_G$ by means of right-invariant
vector fields.

\end{thm}

\section{Proof of the main theorem}

\ssec{The formula}   \label{formula}

Recall (cf. \secref{notation}) that for a right D-module $\M$ on $X$ we have the De Rham
projection $(h\boxtimes \on{id}):\Delta_!(\M)\to \M$.  Let $d_X$ denote the usual De Rham differential
$\M^l\to\M$.

\medskip

Consider the composition
$$\O_G\to J(G)^l\overset{d_X}\to J(G).$$

By construction, $J(G)$ is an $\O_G$-module and the above map is a derivation.
Hence, we obtain a canonical $\O_G$-linear map $\eta:\Omega^1_G\to J(G)$.

For an element $v\in\grg$ let us denote by $v^l$ (resp., $v^r$) the corresponding left-invariant
(resp., right-invariant) vector field on $G$. (Note that the value of $v^l$ (resp., $v^r$) at $1\in G$
is $-v$ (resp., $v$).)
Given a quadratic form $Q$ we will denote by $\eta_Q: \grg\to J(G)$ the composition
$$\grg\overset{v\mapsto v^r}\longrightarrow T_G\overset{Q}\to \Omega^1_G\overset{\eta}\to J(G).$$

\bigskip

The sought-for embedding $\remb:\widetilde{L}_{\grg,Q'}\to \grD_{G,Q}$ is the sum of two D-module maps, 
$\remb_1$ and $\remb_2$. To define $\remb_1$, consider the $\O_{X\times X}$-module map 
$\grg\otimes (\O_X\boxtimes \Omega_X(\Delta))\to \Delta_!(\grD_{G,Q})$ defined as follows:

Let $v\in \grg$ be an element and let $f(x,y)\cdot \omega(y)\in \O_X\boxtimes \Omega_X(\Delta)$ be a local section. Let
us express the right-invariant vector field $v^r$ as
$$v^r=\underset{i}\Sigma\, g_i\cdot v^l_i,\,\,\, g_i\in \O_G.$$ 

We set the image of $v\otimes f(x,y)\cdot \omega(y)$ in $\Delta_!(\grD_{G,Q})$ to be 
$$\underset{i}\Sigma\, \{\,\cdot,\cdot\}(f(x,y)\cdot \lemb(v_i)\boxtimes \omega(y)\cdot g_i),$$
where $\omega\cdot g_i$ is understood as a section of $J(G)\subset \grD^1_{G,Q}$.

It is easy to see that this formula yields $0$ if $f(x,y)\in \O_X\boxtimes \O_X(-\Delta)$. Therefore,
by composing with $(h\boxtimes\on{id})$ we obtain an $\O_X$-module map
$\grg\otimes\F\to \grD_{G,Q}$ and hence a D-module map
$$\grg\otimes(\F\otimes D_X)\to \grD_{G,Q}.$$

To show that it induces a well-defined map of D-modules $\remb_1:\tL'_{\grg,Q'}\to \grD_{G,Q}$, we must prove that
the composition
$\grg\otimes \Omega_X\hookrightarrow \grg\otimes \F\to \Delta_!(\grD_{G,Q})$
factors as
$$\grg\otimes\Omega_X\overset{\rho}\to\Omega_X\hookrightarrow J(G)\hookrightarrow \grD_{G,Q}\to
\Delta_!(\grD_{G,Q}).$$

By property 3(c) of CADO, the above composition sends a section $g\otimes\omega$ to
$$\underset{i}\, \Sigma\,\{\,\cdot,\cdot\}(\lemb(v_i)\boxtimes g_i\cdot \omega)=(\underset{i}\Sigma\, 
\on{Lie}_{v^l_i}(g_i))\cdot \omega.$$

We have the following simple lemma:

\begin{lem}  \label{rho}
For $v\in\grg$ and $v^r=\underset{i}\Sigma\, g_i\otimes v^l_i$ as above, the function
$\underset{i}\Sigma\, \on{Lie}_{v^l_i}(g_i)\in \O_G$ is the constant function equal to $\rho(v)$.
\end{lem}

Thus, we have defined the map $\remb_1$. The map $\remb_2$ is a ``correction term''.  It factors as
$$\tL'_{\grg,Q'}\to \grg\otimes D_X\to J(G) \hookrightarrow \grD_{G,Q},$$
where the middle arrow corresponds to $\eta_{Q_0-Q}:\grg\to J(G)$.

To prove the theorem, we must verify properties (a), (b) and (c), and most importantly, that
$\remb_1+\remb_2$ commutes with the Lie-* bracket. Note, however, that properties (a) and (c) 
are immediate from the construction.

\ssec{Proof of property (b)}   \label{proofofb}

For $u,v\in\grg$, let $\G^Q_{u,v}\in \O_G$ be the function defined by $g\mapsto Q(\on{Ad}_g(u),v)$. Our starting
point is the following observation

\begin{lem}  \label{Killing}
Let $v^r=\underset{i}\Sigma\, g_i\cdot v^l_i,\,\,\, g_i\in \O_G$ as above. Then
$$\on{Lie}_{[u,v_i]^l}(g_i)=\G^{Q_0}_{u,v}.$$
\end{lem}

Fix $f(\cdot,\cdot)$ and $\omega\in\Omega_X$ so that $f(x,y)\cdot \omega(y)\in \O_X\boxtimes\Omega_X(\Delta)$
projects to $1\in\O_X\simeq \O_X\boxtimes\Omega_X(\Delta)/\O_X\boxtimes\Omega_X$. We will use this choice to
have a well-defined $\remb_1(v)\in \grD_{G,Q}$. 

To prove that $\lemb$ and $\remb$ *-commute, it is enough to show that:

\begin{align}
& \{\,\cdot,\cdot\}(\lemb(u)\boxtimes \remb_1(v))=d_1(\G^{Q-Q_0}_{u,v})  \\
& \{\,\cdot,\cdot\}(\lemb(u)\boxtimes \remb_2(v))=d_1(\G^{Q_0-Q}_{u,v}),
\end{align}
where we have used the natural map defined for any D-module $\M$:
$$d_1:\M^l\to \Delta_!(\M),$$
where $\M^l$ is the corresponding left D-module and $d_1$ is the De Rham differential
along the first variable. In our case $\M=J(G)$ we compose it with the embedding $\O_G\to J(G)^l$.

\bigskip

Equation (2) follows from the next result:

\begin{lem} \label{actiononeta}
For $u\in T_G$ and $\omega\in \Omega^1_G$, the Lie-* action $\Theta_G\boxtimes J(G)\to \Delta_!(J(G))$ yields
$$\{\,\cdot,\cdot\}(u\boxtimes \eta(\omega))=\eta(\on{Lie}_u(\omega))-d_1(\langle u,\omega\rangle),$$
where the first term belongs to $J(G)\subset \Delta_!(J(G))$ and
where $\langle u,\omega\rangle\in \O_G$ is the contraction of $u$ and $\omega$.
\end{lem}

This follows immediately from the definitions. Now let us prove Equation (1). 
Consider the expression
\begin{equation} \label{star}
\underset{i}\Sigma\, \{\,\cdot\{\,\cdot,\cdot\}\}(f(y,z)(\lemb(u)
\boxtimes \lemb(v_i)\boxtimes \omega(z)\cdot g_i))\in 
\Delta_!(\grD_{G,Q})  
\end{equation}
on $X\times X\times X$. 
The LHS of Equation (1) is by definition the De Rham under $\on{id}\boxtimes (h\boxtimes\on{id})$ of $(*)$.
  
By applying the Jacobi identity to \eqref{star}, we obtain that the LHS of Equation (1) can be rewritten 
as a sum of three terms:
\begin{align*}
&(a)\;\;\;\;\;\;\;\;\;\;\;\; -(\on{id}\boxtimes h) (\underset{i}\Sigma\, \{\,\cdot,\cdot\} (f(y,x)\cdot\omega(x)\cdot
\on{Lie}_{u^l}(g_i)\boxtimes \lemb(v_i))) \in\grD_{G,Q}\subset\Delta_!(\grD_{G,Q})\\
&(b)\;\;\;\;\;\;\;\;\;\;\;\; \underset{i}\Sigma\, \{\,\cdot,\cdot\}(f(x,y)\cdot\lemb([u,v_i])\boxtimes 
\omega(y)\cdot g_i) 
\in\Delta_!(\grD_{G,Q}) \\
&(c)\;\;\;\;\;\;\;\;\;\;\;\; \underset{i}\Sigma\, Q(u,v_i) 
(\on{id}\boxtimes (h\boxtimes \on{id}))(\on{can}                                            
(f(y,z)\cdot 1'_{x,y}\boxtimes \omega(z)\cdot g_i))\in\Delta_!(\grD_{G,Q}),
\end{align*}
where the notation in term (c) will be explained in \lemref{diagcalc}.

Let us analyze this expression term-by-term. Term (a), which is scheme-theoretically supported
on the diagonal $X\subset X\times X$ equals
$$(\on{id}\boxtimes h) (\underset{i}\Sigma\, \{\,\cdot,\cdot\} 
(f(y,x)\cdot\omega(x)\cdot g_i\boxtimes \lemb([u,v_i]))) \in\grD_{G,Q},$$
because right- and left-invariant vector field commute.

\smallskip

We have the following general assertion:

\begin{prop}  \label{funnybracket}
Let $a,b$ be two sections of a chiral algebra $\A$ such that $\{\,\cdot,\cdot\}(a\boxtimes b)\in \Delta_!(\A)$ is 
scheme-theoretically supported on the diagonal. Let $f(x,y)$ be a local section of $\O_X\boxtimes \O_X(\Delta)$ and let
$\xi$ be the vector field on $X$ equal to the image of $f(x,y)$ under the natural projection 
$\O_X\boxtimes \O_X(\Delta)\to T_X$. Then
$$\{\,\cdot,\cdot\}(f(x,y)(a\boxtimes b))+
(\on{id}\boxtimes h)(\{\,\cdot,\cdot\}(f(y,x)(b\boxtimes a)))=\{\,\cdot,\cdot\}(a\boxtimes b)(\xi,0),$$
where $(\xi,0)$ is the corresponding vector field on $X\times X$.
\end{prop}

\begin{proof}

It suffices to show that $\{\,\cdot,\cdot\}(f(x,y)(a\boxtimes b))-\{\,\cdot,\cdot\}(a\boxtimes b)(\xi,0)\in\Delta_!(\A)$ 
is supported scheme-theoretically on the diagonal and verify that its De Rham projection under $h\boxtimes \on{id}$ 
equals $-(\on{id}\boxtimes h)(\{\,\cdot,\cdot\}(f(y,x)(b\boxtimes a)))$. The latter fact is obvious from the 
anti-symmetry property of the chiral bracket.

Let us multiply $\{\,\cdot,\cdot\}(f(x,y)(a\boxtimes b))-\{\,\cdot,\cdot\}(a\boxtimes b)(\xi,0)$ by a function on 
$X\times X$ of the form $\phi(x)-\phi(y)$. We obtain
$$\{\,\cdot,\cdot\}(f(x,y)\cdot (\phi(x)-\phi(y))(a\boxtimes b))-
\on{Lie}_{\xi}(\phi)(x)\cdot \{\,\cdot,\cdot\}(a\boxtimes b).$$
However, $f(x,y)\cdot (\phi(x)-\phi(y))=\on{Lie}_{\xi}(\phi)\,\on{mod}\, \O_X\boxtimes \O_X(-\Delta)$,
which implies our assertion.

\end{proof}

By applying this proposition to $\lemb([u,v_i])\boxtimes \omega(y)\cdot g_i$, we obtain that
the sum of terms (a) and (b) equals
$$(\underset{i}\Sigma\, \omega\cdot \on{Lie}_{[u,v_i]^l}(g_i))(\xi,0)=-d_1(\on{Lie}_{[u,v_i]^l}(g_i)),$$
which, according to \lemref{Killing}, equals $-d_1(\G^{Q_0}_{u,v})$.

It remains to analyze term (c). We need to show that it equals $d_1(\G^{Q}_{u,v})$.
This results from the following general assertion:

\smallskip
Let $\M$ be a right D-module on $X$ and let $\M^l$ be the corresponding left D-module.

Consider the canonical map of D-modules on $X\times X\times X$
$$\on{can}:(\Delta_2\boxtimes\on{id})_!(j_2{}_*j_2^*(\Omega_X\boxtimes \M))\to \Delta_3{}_!(\M).$$ 
Consider the following section in the RHS:
$$f(y,z)\cdot 1'_{x,y}\boxtimes \omega(z)\cdot m,$$
where $m\in \M^l$, $f(\cdot,\cdot)$ and $\omega$ are as above, and
$1'_{x,y}$ is the section $1'$ defined in \secref{KacMoody} in the first two variables.
 
\begin{lem}  \label{diagcalc}
The $\on{id}\boxtimes (h\boxtimes \on{id})$ projection
of $\on{can}(f(y,z)\cdot 1'_{x,y}\boxtimes \omega(z)\cdot m)$ equals $-d_1(m)$.
\end{lem}

\begin{proof}

Let us multiply both sides by a function on $X\times X$ of the form $\phi(x)-\phi(y)$. In both cases
we get $d\phi\cdot m\in M\subset \Delta_!(M)$. Therefore, the expression
\begin{equation} \label{2star}
(\on{id}\boxtimes (h\boxtimes \on{id}))
(\on{can}(f(y,z)\cdot 1'_{x,y}\boxtimes \omega(z)\cdot m))+d_1(m)\in\Delta_!(\M)
\end{equation}
is killed by the equation of the diagonal.
However, since
$$(h\boxtimes \on{id})\circ (\on{id}\boxtimes(h\boxtimes \on{id}))(\on{can}(f(y,z)\cdot 1'_{x,y}\boxtimes \omega(z)\cdot
m))=0= (h\boxtimes \on{id})(-d_1(m)),$$
we obtain that \eqref{2star} vanishes identically.

\end{proof}

\ssec{Proof of the compatibility with the Lie-* bracket}     \label{proofofcom}

Let $v,v'\in\grg$ be two elements. The expression $\{\,\cdot,\cdot\}(\remb(v)\boxtimes \remb(v'))$ is a well-defined
section of $\Delta_!(\grD_{G,Q})$ and $\{\,\cdot,\cdot\}(v\boxtimes v')$ is a well-defined section of
$\Delta_!(\tL'_{\grg,Q'})$. Our goal is to prove that
\begin{equation} \label{3star}
\{\,\cdot,\cdot\}(\remb(v)\boxtimes \remb(v'))-\remb(\{\,\cdot,\cdot\}(v\boxtimes v'))=0\in \Delta_!(\grD_{G,Q})
\end{equation}

First, we will prove that the RHS of \eqref{3star} is killed by the equation of the diagonal. We must show that
for a function $\phi\in\O_X$,
$$(\phi(x)-\phi(y))\cdot \{\,\cdot,\cdot\}(\remb(v)\boxtimes \remb(v'))=d_X(\phi)\cdot (-Q_0-Q)(v,v').$$

It is clear that $\{\,\cdot,\cdot\}(\remb_2(v)\boxtimes \remb_2(v'))$ equals zero. We can compute
$\{\,\cdot,\cdot\}(\remb_1(v)\boxtimes \remb_2(v'))$ using \lemref{actiononeta}. We obtain:
\begin{align*}
&(\phi(x)-\phi(y))\cdot \{\,\cdot,\cdot\}(\remb_1(v)\boxtimes \remb_2(v'))=(\phi(x)-\phi(y))\cdot
d_1((Q-Q_0)(v,v'))= \\
&d_X(\phi)\cdot (Q_0-Q)(v,v').
\end{align*}
Similarly, $(\phi(x)-\phi(y))\cdot \{\,\cdot,\cdot\}(\remb_2(v)\boxtimes \remb_1(v'))=
d_X(\phi)\cdot (Q_0-Q)(v,v')$. Thus, it remains to compute $\{\,\cdot,\cdot\}(\remb_1(v)\boxtimes \remb_1(v'))$.

The latter is the $(h\boxtimes \on{id})\boxtimes (h\boxtimes \on{id})$ projection of
\begin{equation} \label{4star}
\underset{i,j}\Sigma\,\{\{\,\cdot,\cdot \},\{\,\cdot,\cdot\}\}
(f(x,y)\cdot \lemb(v_i)\boxtimes \omega(y)\cdot g_i\boxtimes f(z,w)\cdot \lemb(v'_j)\boxtimes \omega(w)\cdot g'_j)
\end{equation}
where $f(\cdot,\cdot)\in \O_X\boxtimes\O_X(\Delta)$ and $\omega\in\Omega_X$ are as before.

\smallskip

Using Jacobi identity, we can rewrite \eqref{4star} as a sum over $i,j$ of four terms:
\begin{align*}
&(a)\;\;\;\;\;\;\;\;\;\; \{\,\cdot,\{\{\,\cdot,\cdot\},\cdot\}\}(f(x,y)\cdot 
\lemb(v_i)\boxtimes \omega(y)\cdot g_i\boxtimes
f(z,w)\cdot \lemb(v'_j)\boxtimes \omega(w)\cdot g'_j) \\
&(b)\;\;\;\;\;\;\;\;\;\; -\sigma_{3,4} \{\,\cdot,\{\{\,\cdot,\cdot\},\cdot\}\}(f(x,y)\cdot \lemb(v_i)\boxtimes 
\omega(y)\cdot g_i\boxtimes f(w,z)\cdot\omega(z)\cdot g'_j\boxtimes \lemb(v'_j)) \\
&(c)\;\;\;\;\;\;\;\;\;\; -\sigma_{1,2} \{\,\cdot,\{\{\,\cdot,\cdot\},\cdot\}\}(f(y,x)\cdot \omega(x)\cdot
g_i\boxtimes \lemb(v_i) \boxtimes f(z,w)\cdot \lemb(v'_j)\boxtimes \omega(w)\cdot g'_j)  \\
&(d)\;\;\;\;\;\;\;\;\;\; \sigma_{1,2} \circ\sigma_{3,4}
\{\,\cdot,\{\{\,\cdot,\cdot\},\cdot\}\}(f(y,x)\cdot \omega(x)\cdot g_i\boxtimes \lemb(v_i)\boxtimes f(w,z)\cdot
\omega(z)\cdot g'_j\boxtimes \lemb(v'_j)),
\end{align*}
where $\sigma_{1,2},\sigma_{3,4}\in S^4$ are the transpositions of the corresponding factors.

It is easy to see that term (b) vanishes identically and that term (a) is killed by $(\phi(y)-\phi(w))$.
The contribution of term (c), multiplied by $(\phi(y)-\phi(w))$ consists of two parts. 
One is $\on{id}\boxtimes (h\boxtimes \on{id})$ applied to
$$-\underset{i,j}\Sigma\, \{\,\cdot,\{\,\cdot,\cdot\}\}(f(y,x)\cdot 
(\phi(x)-\phi(z))\cdot \omega(x)\cdot g_i\boxtimes \lemb([v_i,v'_j])\boxtimes f(y,z)\cdot \omega(z)\cdot g'_j),$$
and the other is $(\on{id}\boxtimes h)\boxtimes (h\boxtimes \on{id})$ applied to
$$-\underset{i,j}\Sigma\, Q(v_i,v'_j)\cdot  \{\,\cdot,\cdot\}
(f(y,x)\cdot (\phi(x)-\phi(w))\cdot \omega(x)\cdot g_i\boxtimes \on{can}(f(z,w)\cdot 1'_{y,z}\boxtimes 
\omega(w)\cdot g'_j)),$$
where "$\on{can}$" is as in \lemref{diagcalc}.

The first part can be computed using \lemref{Killing} and we obtain
\begin{align*}
&-\underset{i,j}\Sigma\, \{\,\cdot,\{\,\cdot,\cdot\}\}(f(y,x)\cdot 
(\phi(x)-\phi(y))\cdot \omega(x)\cdot g_i\boxtimes \lemb([v_i,v'_j])\boxtimes f(y,z)\cdot \omega(z)\cdot g'_j) \\
&-\underset{i,j}\Sigma\, \{\,\cdot,\{\,\cdot,\cdot\}\}(f(y,x)\cdot 
(\phi(y)-\phi(z))\cdot \omega(x)\cdot g_i\boxtimes \lemb([v_i,v'_j])\boxtimes f(y,z)\cdot \omega(z)\cdot g'_j)=\\
&=\underset{i,j}\Sigma\, d_X(\phi)\cdot g_i\cdot \on{Lie}_{[v_i,v'_j]^l}(g'_j)+
\underset{i,j}\Sigma\, d_X(\phi)\cdot \on{Lie}_{[v'_j,v_i]^l}(g_i)\cdot g'_j=\\
&=d_X(\phi)\cdot (\underset{i}\Sigma\, g_i\cdot \G^{Q_0}_{v_i,v'}+\underset{j}\Sigma\, g'_j\cdot \G^{Q_0}_{v'_j,v})
=-2d_X(\phi)\cdot Q_0(v,v').
\end{align*}

The second part of term (c) can be calculated using \lemref{diagcalc} and we obtain:
$$\underset{i,j}\Sigma\, d_X(\phi) \cdot g_i\cdot g'_j\cdot Q(v_i,v'_j)=d_X(\phi)\cdot Q(v,v').$$

Finally, term (d) multiplied by $(\phi(y)-\phi(w))$ gives
\begin{align*}
&-\underset{i,j}\Sigma\, \{\,\cdot,\{\,\cdot,\cdot\}\}(d_X(\phi)(x)\cdot g_i\boxtimes f(z,y)\cdot 
\omega(y)\cdot \on{Lie}_{v_i^l}(g'_j)\boxtimes \lemb(v'_j))= \\
&=-d_X(\phi)\cdot \underset{i,j}\Sigma\, \on{Lie}_{v_i^l}(g_j')\cdot \on{Lie}_{v'_j{}^l}(g_i)=
d_X(\phi)\cdot \underset{i,j}\Sigma\, g_i\cdot \on{Lie}_{[v_i,v'_j]^l}(g'_j)=-d_X(\phi)\cdot Q_0(v,v').
\end{align*}

\smallskip

By summing up, we obtain the desired
$$(\phi(x)-\phi(y))\cdot \{\,\cdot,\cdot\}(\remb_1(v)\boxtimes \remb_1(v'))=d_X(\phi)\cdot (-Q-Q_0)(v,v').$$

\medskip

Thus, the LHS of \eqref{3star} can be regarded as a map
$$\psi:\grg\otimes\grg\to \grD_{G,Q}.$$
By property 3 of CADO, the image of $\psi$ is contained in $\grD_{G,Q}^0=J(G)$. In addition, according to
property (b) of the map $\remb$ that has already been checked, the image of $\psi$ *-commutes with 
$\lemb:\widetilde{L}_{\grg,Q}\to \grD_{G,Q}$.

However, since the space of $\grg\otimes \hat\O_x$-invariants in $\O_{G(\hat\O_x)}$ consists of constant functions,
a section of $J(G)$   *-commutes with the image of $\lemb$ if and only if it belongs to $\Omega_X\subset J(G)$. 
Hence, $\psi$ is a map $\grg\otimes\grg\to\Omega_X$.

To prove that it vanishes we will use a "conformal dimension" argument:

\smallskip

Since our constructions are local, we can assume that $X=\CC$,
and $\psi$ is invariant with respect to the group of all automorphisms of $\CC$.
Hence, the image of $\psi$ belongs belongs to the space of $\on{Aut}(\CC)$-invariant sections of
$\Omega_{\CC}$, but there are none.

\section{Relation to semi-infinite cohomology}   \label{sisec}

\ssec{The BRST complex}

Our goal in this section is to compute the semi-infinite cohomology 
of ${\mathbb V}_{G,Q}=(\grD_{G,Q})_x$ viewed as a
module over the affine algebra. First, let us recall the basic definitions 
concerning semi-infinite cohomology in
the context of chiral algebras, following \cite{BD}, Section 3.8.

Let $\M$ be a finitely generated locally free right D-module and 
let $$\M^*:=\Hom(\M,D_X\otimes\Omega_X)$$ be its
Verdier dual. The direct sum
$\Omega_X\oplus (\M\oplus\M^*)$
has a natural structure of a (super) Lie-* algebra, where the Lie-* bracket 
vanishes on $\Omega_X$ and
$\M\boxtimes \M^*\to\Delta_!(\Omega_X)$
is the canonical map, cf. \cite{BD}, \cite{Ga}.

\smallskip  

We define the Clifford (super) chiral algebra $\on{Cliff}(\M)$ as the 
quotient of the universal
enveloping chiral algebra 
$\U(\Omega_X\oplus (\M\oplus\M^*))$
by the ideal generated by the anti-diagonal copy of $\Omega_X$. The 
$\ZZ_2$-grading on $\on{Cliff}(\M)$
can be extended to a $\ZZ$-grading, by letting $\M$ 
(resp., $\M^*$, $\Omega_X$) have degree $-1$ (resp., $1$, $0$).

Consider the $0$-th graded component of the second associated graded quotient with respect to the
canonical filtration on $\on{Cliff}(\M)$, $\on{gr}_2(\on{Cliff}(\M))^0$. It has a natural structure 
of a Lie-* algebra isomorphic to $\on{End}(\M):=\M\overset{!}\otimes \M^*$, cf. {\it loc. cit.}.
Hence, if we do not mod out by $\Omega_X\simeq \on{Cliff}(\M)_1^0$, we obtain an extension of Lie-* algebras
$$0\to \Omega_X \to \widetilde{\on{End}}(\M) \to \on{End}(\M)\to 0.$$

\medskip

Now let us assume that $\M=L$ is a Lie-* algebra, which is locally free as a D-module. As is explained
in \cite{BD} or \cite{Ga}, the Lie-* bracket on $L$ yields a Lie-* algebra map
$L\to \on{End}(L)$. In particular, from $\widetilde{\on{End}}(L)$ we obtain a canonical central extension
$$0\to\Omega_X\to \widetilde{L}_{can}\to L\to 0.$$
Let us denote by $\widetilde{L}_{-can}$ the Baer negative of the extension $\widetilde{L}_{can}$. We will
denote by $\widetilde{\U}(L)_{can}$ (resp., $\widetilde{\U}(L)_{-can}$) the quotient of 
$\U(\widetilde{L}_{can})$ (resp., $\U(\widetilde{L}_{-can})$) by the antidiagonal copy of $\Omega_X$.

\smallskip

Let $\N_x$ be a representation of the extension
$$0\to \CC\to DR(\D^*_x,\widetilde{L}_{-can})\to DR(\D_x^*,L)\to 0.$$
According to \lemref{catmodules},, such $\N_x$ is a vector space underlying 
a chiral module over $\widetilde{\U}(L)_{-can}$. 

It is explained in {\it loc. cit.} that the tensor product $\N_x\otimes \on{Cliff}(L)_x$,
which is a $\ZZ$-graded vector space, acquires a differential $\delta$ (called BRST) of degree $1$.

The semi-infinite cohomology of $\N_x$ with respect to $L$ is defined as the cohomology of this complex:
$$\Hsik(L,\N_x):=\on{H}^k(\N_x\otimes \on{Cliff}(L)_x,\delta).$$ 

\ssec{Identification of the canonical extension}

From now on, we will take $L=L_{\grg}=\grg\otimes D_X$. According to the previous subsection, there exists
a canonical central extension of Lie-* algebras
$$0\to \Omega_X\to \widetilde{L_{\grg}}_{-can} \to L_{\grg}\to 0.$$

Recall also that in \secref{mainresult} we defined another central extension of $L_{\grg}$ corresponding
to an invariant form $Q$, namely $\tL'_{\grg,Q}$. 
In this subsection our goal will be to prove the following assertion:

\begin{thm}   \label{extension}
There is a canonical isomorphism of extensions $\tL'_{\grg,-Q_0}\simeq 
\widetilde{L_{\grg}}_{-can}$.
\end{thm}

\begin{proof}

For $L=L_{\grg}$, $\on{End}(L_{\grg})\simeq \on{End}(\grg)\otimes D_X$ and the Lie-* bracket
on it comes from the usual Lie algebra structure on $\on{End}(\grg)$. The map
$L_{\grg}\to \on{End}(L_{\grg})$ comes from the action map $\grg\to \on{End}(\grg)$.

\smallskip

To prove the theorem, it suffices
to construct a Lie-* algebra map
$$\temb:\tL'_{\grg,-Q_0}\to \widetilde{\on{End}}(\L),$$
such that the induced map from $\Omega_X\subset \tL'_{\grg,-Q_0}$ to 
$\Omega_X\subset \widetilde{\on{End}}(\L)\subset \on{Cliff}(L)$ is $-\on{id}$.

\medskip

let $v\in \grg$ and let $v^r=\underset{i}\Sigma\, g_i\cdot v^l_i,\,\,\, g_i\in \O_G$. Then the image
of $v$ in $\on{End}(\grg)\simeq \grg\otimes \grg^*$ is the value at $1\in G$ of the element
$$-\underset{i}\Sigma\, v_i\otimes d_G(g_i) \in \grg\otimes \Omega^1_G.$$

The required map $\temb$ is defined by a formula similar to the one of \secref{formula}. Namely it sends
a section $v\otimes f(x,y)\cdot \omega(y)\in \grg\otimes \F\otimes D_X$ to the $(h\boxtimes \on{id})$
projection of 
$$\underset{i}\Sigma\, \{\,\cdot,\cdot\}(f(x,y)\cdot\omega(y)\cdot v_i\boxtimes d_G(g_i)(1))\in 
\Delta_!(\on{Cliff}(L_{\grg})).$$

The fact that $\temb|_{\Omega_X}$ is $-\on{id}$ follows immediately from \lemref{rho}. To show that
it commutes with the Lie-* algebra bracket, we proceed as in the proof of \thmref{main}. Namely,
by the same argument as in \secref{proofofcom}, it suffices to show that the difference

\begin{equation} \label{nstar}
\temb(\{\,\cdot,\cdot\}(v\boxtimes v'))-\{\,\cdot,\cdot\}(\temb(v)\boxtimes \temb(v'))\in 
\Delta_!(\on{Cliff}(L_{\grg}))
\end{equation}
is killed by a function of the form $\phi(x)-\phi(y)\in \O_{X\times X}$.

\smallskip

By applying the (super-) Jacobi identity, we obtain that the second term in \eqref{nstar} is the projection under 
$(h\boxtimes \on{id})\boxtimes (h\boxtimes \on{id})$ of the sum over $i$ and $j$ of four terms
\begin{align*}
&\{\,\cdot,\{\{\,\cdot,\cdot\},\cdot\}\}(f(x,y)\cdot v_i\boxtimes \omega(y)\cdot d_G(g_i)(1)\boxtimes
f(z,w)\cdot v'_j\boxtimes \omega(w)\cdot d_G(g'_j)(1)) \\
&\sigma_{3,4} \{\,\cdot,\{\{\,\cdot,\cdot\},\cdot\}\}(f(x,y)\cdot v_i\boxtimes \omega(y)\cdot d_G(g_i)(1)\boxtimes
f(w,z)\cdot \omega(z)\cdot d_G(g'_j)(1) \boxtimes v'_j) \\
&\sigma_{1,2} \{\,\cdot,\{\{\,\cdot,\cdot\},\cdot\}\}(f(y,x)\cdot \omega(x)\cdot
d_G(g_i)(1)\boxtimes v_i \boxtimes f(z,w)\cdot v'_j\boxtimes \omega(w)\cdot d_G(g'_j)(1)) \\
&\sigma_{1,2} \circ\sigma_{3,4}
\{\,\cdot,\{\{\,\cdot,\cdot\},\cdot\}\}(f(y,x)\cdot \omega(x)\cdot d_G(g_i)(1)\boxtimes v_i\boxtimes f(w,z)\cdot
\omega(z)\cdot d_G(g'_j)(1)\boxtimes v'_j).
\end{align*}

Let us analyze the above expression term-by-term. First of all, 
it is easy to see that the second and the third terms vanish. Secondly, the first term is 
killed by multiplication by $\phi(x)-\phi(y)$. The fourth term yields that:
\begin{align*}
&(\phi(x)-\phi(y))\cdot\{\,\cdot,\cdot\}(\temb(v)\boxtimes \temb(v'))=-d_X(\phi)\cdot \underset{i,j}\Sigma\,
(\on{Lie}_{v_i^l}(g'_j)\cdot \on{Lie}_{v'_j{}^l}(g_i))(1)=\\
&-d_X(\phi)\cdot \underset{i,j}\Sigma\, g'_j\cdot \on{Lie}_{[v'_j,v_i]^l}(g_i)
=-d_X(\phi)\cdot \underset{i}\Sigma\, g'_j\cdot \G^{Q_0}_{v'_j,v}=
d_X(\phi)\cdot Q_0(v,v').
\end{align*}

However, $(\phi(x)-\phi(y))\cdot \temb(\{\,\cdot,\cdot\}(v\boxtimes v'))=d_X(\phi)\cdot Q_0(v,v')$,
which is what we had to prove.

\end{proof}

\ssec{Computation of $\Hsi(\grD_{G,Q})$}

In what follows, we will denote by $\tg'_Q$ the affine algebra $DR(\D_x^*,\tL'_{\grg,Q})$. (Non-canonically,
$\tg'_Q$ and $\tg_Q$ are isomorphic, but this isomorphism does not respect the action of the group 
$\on{Aut}(\D_x)$ of automorphisms of the formal disc, which acts on the whole picture.) In particular,
if $\M^1_x$ is a $\tg_{Q_1}$-representation and $\M^2_x$ is a $\tg'_{Q_2}$-representation, the tensor 
product $\M^1_x\otimes \M^2_x$ is a $\tg'_{Q_1+Q_2}$-representation via the diagonal action.

Let $\M_x$ be a (continuous) module over the affine algebra $\tg_Q$ and let us 
consider the tensor product $\M_x\otimes {\mathbb V}_{G,Q}$. We let $\widetilde{\grg}_Q$ act on it via
the embedding $\lemb:L_{\grg,Q}\to \grD_{G,Q}$. However, since by \thmref{main}, ${\mathbb V}_{G,Q}$
is a bi-module over the affine algebra, the above tensor product has an additional structure. 
Namely, by combining \thmref{main} and \thmref{extension}, we obtain on it an action of $\tg'_{-Q_0}$.

Hence, it makes sense to consider $\Hsi(L_{\grg},\M_x\otimes {\mathbb V}_{G,Q})$, and will carry 
a continuous $\tg_Q$-action.

\begin{thm} \label{semiinf}
Assume that $\M_x$ is such that the action of $\grg\otimes\hat\O_x\subset \tg_Q$ on it
integrates to an $G(\hat\O_x)$-action. Then there is a canonical isomorphism of $\tg_Q$-modules:
$$\Hsik(L_{\grg},\M_x\otimes {\mathbb V}_{G,Q})\simeq \M_x\otimes \on{H}^k(\grg,\CC).$$
\end{thm}

The proof will consist of two steps. First, we will consider the case $Q=0$ and $\M_x=\CC$.

\begin{prop}  \label{semivacuum}
There is a canonical isomorphism of $\grg\otimes\hat\K_x$-modules 
$$\Hsik(L_{\grg},{\mathbb V}_{G,Q})\simeq \on{H}^k(\grg,\CC),$$ where the action of $\grg\otimes\hat\K_x$
on the RHS is trivial.
\end{prop}

\begin{proof}

We will use the following statement, valid for an arbitrary locally free Lie-* algebra:

Let $\N_x$ be a module over $DR(\D^*_x,\widetilde{L}_{-can})$, which is induced from a
$DR(\D_x,L)$-module, i.e.
$$\N_x\simeq \on{Ind}^{DR(\D^*_x,\widetilde{L}_{-can})}_{DR(\D_x,L)}(\overline{N}),$$
where induction is understood in the restricted sense, i.e. $1\in\CC\subset DR(\D_x^*,\widetilde{L}_{-can})$
acts as identity.

\begin{lem}
Under the above circumstances there is a canonical isomorphism
$$\Hsik(L,\N_x)\simeq \on{H}^k(DR(\D_x,L),\overline{N}).$$
\end{lem}

We will apply this lemma in the case $\N_x\simeq {\mathbb V}_{G,0}$. It is easy to see, as in
\lemref{semidual}, that the embedding $J(G)\to \grD_{G,0}$ induces an isomorphism of $\tg'_{-Q_0}$-modules
$${\mathbb V}_{G,0}\simeq \on{Ind}^{\tg'_{-Q_0}}_{\grg\otimes\hat\O_x}(\O_{G(\hat\O_x)}).$$

Moreover, the restriction of the $\grg\otimes\hat\K_x$-action on the LHS to
$\grg\otimes\hat\O_x$, that comes from the embedding $\lemb:L_{\grg}\to \grD_{G,0}$ coincides
with the action that comes from the natural $\grg\otimes\hat\O_x$-action on $\O_{G(\hat\O_x)}$ by
left-invariant vector fields. This follows from the fact that the embeddings $\lemb$ and $\remb$  *-commute with
one another.

Hence, we obtain an isomorphism
$\Hsik(L,{\mathbb V}_{G,0})\simeq \on{H}^k(\grg\otimes\hat\O_x,\O_{G(\hat\O_x)})$
and it respects the $\grg\otimes\hat\O_x$-action on both sides. Since the kernel of the evaluation map
$G(\hat\O_x)\to G$ is pro-unipotent, we have an isomorphism of $\grg\otimes\hat\O_x$-modules
$$\on{H}^k(\grg\otimes\hat\O_x,\O_{G(\hat\O_x)})\simeq \on{H}^k(\grg,\O_G)\simeq \on{H}^k(\grg,\CC),$$
where the $\grg\otimes\hat\O_x$-action on the RHS is trivial.

\medskip

To prove the proposition, it remains to show that the $\grg\otimes\hat\K_x$-action on $\on{H}^k(\grg,\CC)$
is trivial as well. Consider the group $\on{Aut}(\D_x)$ of automorphisms of the formal disc. This
group acts on the whole picture. In particular, our homomorphism
$$\grg\otimes\hat\K_x\to \on{End}(\on{H}^k(\grg,\CC))$$
is $\on{Aut}(\D_x)$-equivariant. However, it is easy to see that any such homomorphism, which is, moreover,
trivial on $\grg\otimes\hat\O_x\subset \grg\otimes\hat\K_x$ is zero.

\end{proof}

The second step in the proof of \thmref{semiinf} is the following result:

\medskip

Let $\M_x$ be a $\tg_Q$-module as above and let us consider the tensor product
$\M_x\otimes {\mathbb V}_{G,Q'}$, where $Q'$ is another invariant form.
It has a natural bi-module structure with respect to $\tg_{Q'}$ and $\tg'_{Q-Q_0-Q'}$:

We let $\tg_{Q'}$ act via $\lemb$ on ${\mathbb V}_{G,Q'}$, and 
$\tg'_{Q-Q_0-Q'}$ will act diagonally via $\remb$ and the $\tg_Q$-action on $\M_x$.

Consider now the tensor product $\M_x\otimes {\mathbb V}_{G,Q'-Q}$. It has a bi-module structure with
respect to the same Lie algebras, where $\tg'_{Q'}$ acts diagonally (via $\lemb$) and
$\tg'_{Q-Q_0-Q'}$ acts only on ${\mathbb V}_{G,Q'-Q}$ (via $\remb$).

\begin{thm}  \label{regular}
Assume that $\M_x$ is such that the action of $\grg\otimes\hat\O_x\subset \tg_Q$ on it
integrates to an $G(\hat\O_x)$-action. Then the above bi-modules $\M_x\otimes {\mathbb V}_{G,Q'}$
and $\M_x\otimes {\mathbb V}_{G,Q'-Q}$ are canonically isomorphic.
\end{thm}

It is clear that \thmref{regular} combined with \propref{semivacuum} imply \thmref{semiinf}, when
one takes $Q'=Q$.

\ssec{Proof of \thmref{regular}}

Recall that if $\A_1$ and $\A_2$ are chiral algebras, then their tensor product
$\A_1\overset{!}\otimes\A_2$ is also naturally a chiral algebra.

Thus, let us consider the chiral algebras $\A_{\grg,Q}\overset{!}\otimes \grD_{G,Q'}$ and 
$\A_{\grg,Q}\overset{!}\otimes \grD_{G,Q'-Q}$. We will denote by $\iemb$ the embedding of the Lie-*
algebra $\tL_{\grg,Q}$ into each of them along the first factor. Note that in addition, one has the
natural (diagonal) maps
$$\lemb+\iemb:\tL_{\grg,Q}\to \A_{\grg,Q}\overset{!}\otimes \grD_{G,Q'-Q} \text{ and }
\remb+\iemb:\tL'_{\grg,Q-Q_0-Q'}\to \A_{\grg,Q}\overset{!}\otimes \grD_{G,Q'}.$$

\begin{prop}
There exists a canonical isomorphism 
$$\phi:\A_{\grg,Q}\overset{!}\otimes \grD_{G,Q'}\to \A_{\grg,Q}\overset{!}\otimes \grD_{G,Q'-Q},$$
such that 

\smallskip

\noindent {\em (a)} $\phi$ commutes with the embedding of $J(G)$ into 
$\grD_{G,Q'}\subset \A_{\grg,Q}\overset{!}\otimes \grD_{G,Q'}$ and
$\grD_{G,Q'-Q}\subset \A_{\grg,Q}\overset{!}\otimes \grD_{G,Q'-Q}$.

\smallskip

\noindent {\em (b)} $\phi\circ \lemb=\lemb+\iemb$.

\smallskip

\noindent {\em (c)} $\phi\circ (\remb+\iemb)=\remb$.

\end{prop}

\begin{proof}

By construction, $\A_{\grg,Q}\overset{!}\otimes \grD_{G,Q'}$ is a quotient of the universal
enveloping chiral algebra of $\tL_{\grg,Q}\oplus J(G)\oplus L_{\grg}$. Therefore, $\phi$, if it exists,
is uniquely described by a Lie-* algebra map
$$\tL_{\grg,Q}\oplus J(G)\oplus L_{\grg}\to \A_{\grg,Q}\overset{!}\otimes \grD_{G,Q'-Q}.$$

On the second and the third summands $\phi$ is determined by conditions (a) and (b), respectively. 
The restriction of $\phi$ to $\Omega_X\subset \tL_{\grg,Q}$ is the identity map onto $\Omega_X$ in
the target and the restriction to $\grg\otimes D_X\subset \tL_{\grg,Q}$ is determined by the
following condition:

The adjoint action gives rise to a map $\grg\to \grg\otimes \O_G$. From it we obtain a D-module map
$$\grg\otimes D_X\to (\grg\otimes D_X)\overset{!}\otimes J(G).$$
(On the level of fibers, this map is the co-action 
$\grg\otimes \hat\K_x/\hat\O_x\to \grg\otimes \hat\K_x/\hat\O_x\otimes \O_{G(\hat\O_x)}$.)
We need that the diagram
$$
\CD
\grg\otimes D_X  @>>> (\grg\otimes D_X)\overset{!}\otimes J(G) \\
@VVV        @VVV     \\
\tL_{\grg,Q}   @>{\phi}>>  \tL_{\grg,Q}\overset{!}\otimes \grD_{G,Q'-Q}
\endCD
$$
commutes.

\smallskip

The fact that the resulting map of D-modules commutes with the Lie-* bracket is straightforward.
Thus, one obtains a map of chiral algebras
$$\U(\tL_{\grg,Q}\oplus J(G)\oplus L_{\grg})\to \A_{\grg,Q}\overset{!}\otimes \grD_{G,Q'-Q},$$
and one easily checks that it factors through the natural surjection
$$\U(\tL_{\grg,Q}\oplus J(G)\oplus L_{\grg})\twoheadrightarrow \A_{\grg,Q}\overset{!}\otimes \grD_{G,Q'}.$$

The fact that condition (c) is satisfied follows from the formula for the embedding $\remb$, cf. \secref{formula}.
The fact that $\phi$ is an isomorphism is easy, since one can construct its inverse by a similar
procedure.

\end{proof}

Now, let $\M_x$ be as in the formulation of the theorem. The tensor products
$\M_x\otimes {\mathbb V}_{G,Q'}$ and $\M_x\otimes {\mathbb V}_{G,Q'-Q}$ are 
chiral modules supported at $x\in X$ for $\A_{\grg,Q}\overset{!}\otimes \grD_{G,Q'}$
and $\A_{\grg,Q}\overset{!}\otimes \grD_{G,Q'-Q}$, respectively.

To prove the theorem, it suffices to construct an isomorphism 
$$\phi_{\M_x}:\M_x\otimes {\mathbb V}_{G,Q'}\to \M_x\otimes {\mathbb V}_{G,Q'-Q},$$
which covers the isomorphism $\phi$ of the above proposition.

\smallskip

However, it is easy to see that there exists a well-defined map $\phi_{\M_x}$ such that its restriction
to $$\M_x\otimes 1_x\subset \M_x\otimes (\grD_{G,Q})_x\simeq \M_x\otimes {\mathbb V}_{G,Q'}$$ is the map
$$\M_x\to \M_x\otimes \O_{G(\hat\O_x)}=\M_x\otimes J(G)_x\subset \M_x\otimes (\grD_{G,Q'-Q})_x=\M_x\otimes {\mathbb
V}_{G,Q'-Q},$$ where the first arrow is the co-action map corresponding to the $G(\hat\O_x)$-action on $\M_x$.

\section{Appendix: From $\grD_{G,Q}$-modules to D-modules}
\label{appendix}

The discussion in this section is substantially based on the
unpublished manuscript \cite{BD1}, Sect. 7, where the formalism
of D-modules on ind-schemes is developed. For that reason,
we call this section an Appendix, and its results and techniques
are not directly connected to the contents of the main part of the paper.

\ssec{A characterization of chiral $\grD_{G,Q}$-modules}

In this section the point $x\in X$ will be fixed and we will
denote $\hat\O_x$ by $\CC[[t]]$ and $\hat\K_x$ by $\CC((t))$,
respectively. Moreover, all chiral $\grD_{G,Q}$-modules we that will
consider are supported at $x$; therefore, by 
a slight abuse of notation, we will identify a chiral module
$\M$ with the corresponding vector space $i_x^!(\M)[1]$.

Recall that $G[[t]]$ has a canonical structure of a group-scheme
and $G((t))$ of a group-indscheme. In particular, 
$\O_{G((t))}$ is a topological commutative algebra, which is
a continuous bi-module over $\grg((t))$.

We have the following proposition:

\begin{prop}   \label{catrep}
To specify a structure of a chiral $\grD_{G,Q}$-module
on a vector space $\M$ is the same as to endow it with 
continuous (w.r. to the discrete topology on $\M$) actions of
$\O_{G((t))}$ and $\tg_Q$ compatible in the sense that for $\xi\in\tg_Q$,
$f\in \O_{G((t))}$ and $m\in \M$,
$$\xi\cdot(f\cdot m)=f\cdot (\xi\cdot m)+\on{Lie}_{\xi^l}(f)\cdot m,$$
where $\xi^l$ is the corresponding left-invariant vector field on $G((t))$.
\end{prop}

\begin{proof}

First, let $\B$ be a commutative chiral algebra. Consider
the projective limit $\hat \B_x:=
\underset{\longleftarrow}{lim}(\B^\ell_i)_x$, 
where $\B_i$ runs over the set of all chiral subalgebras
$\B_i\subset \B$ with $\B_i|_{X-x}\simeq \B|_{X-x}$.

Then $\hat \B_x$ is a commutative algebra, which carries 
a natural topology. Moreover, a structure of 
a chiral $\B$-module on a vector space $\M$ it amounts to a 
continuous $\hat \B_x$-action on $\M$. 

When $\B$ is of the form $\B=J(\C)$, for a q.c. sheaf of 
$\O_X$-algebras $\C$, $\hat \B_x$ represents the functor
on the category of $\CC$-algebras given by $A\to \Hom_{\O_X}(\C,A((t))$. 
In particular, $\widehat{J(G)}_x\simeq \O_{G((t))}$. 

\medskip

Let us suppose now that $\B$ is a Lie-* module over
a Lie-* algebra $L$. Then we obtain a natural continuous action
of $DR(\D^*_x,L)$ on $\hat \B_x$. It is easy to see that in the
above example of $\B=J(G)$ and the action
$\grg\otimes D_X\boxtimes J(G)\to\Delta_!(J(G))$
of \secref{secconstr}, the resulting $\grg((t))$-action
on $\O_{G((t))}$ coincides with the natural action
by left-invariant vector fields.

\medskip

That said, the assertion of the proposition becomes a direct
corollary of the construction of $\grD_{G,Q}$ as a quotient
of the universal enveloping algebra of $J(G)\oplus L_{\grg}$
combined with \lemref{catmodules}.

\end{proof}

Let us observe now that \thmref{main} implies that every 
chiral $\grD_{G,Q}$-module $\M$ 
carries also an action of $\tg'_{Q'}$, which commutes
with the initial $\tg_Q$-action. In what follows, we will refer to the
above $\tg_Q$ and $\tg'_{Q'}$-actions as ``left'' and ``right'', 
respectively. 

\begin{lem}
For $\M$ as above, the right $\tg'_{Q'}$-action is compatible with the
$\O_{G((t))}$-action in the sense that for $\xi'\in\tg_Q$,
$f\in \O_{G((t))}$ and $m\in \M$,
$$\xi\cdot(f\cdot m)=f\cdot (\xi\cdot m)+\on{Lie}_{\xi^r}(f)\cdot m,$$
where $\xi^r$ is the corresponding right-invariant vector field on $G((t))$.
\end{lem}

\begin{proof}

Let $L_{\grg}\to \Theta_G$ be the Lie-* algebra map
corresponding to the embedding $\grg\to T_G$
by means of right-invariant vector fields.
Let is consider the corresponding Lie-* bracket
$$L_{\grg}\boxtimes J(G)\to \Theta_G\boxtimes J(G)
\to\Delta_!(J(G)).$$

As in the proof of \propref{catrep} above, we need to show that
if we compose the above Lie-* bracket with the natural surjection
$\wt{L}'_{\grg,Q'}\to L_{\grg}$, we obtain a commutative
diagram
$$
\CD
\wt{L}'_{\grg,Q'}\boxtimes J(G) @>>>  \Delta_!(J(G)) \\
@V{\remb\boxtimes \on{id}}VV            @VVV        \\
\grD_{G,Q}\boxtimes \grD_{G,Q} @>>> \Delta_!(\grD_{G,Q}).
\endCD
$$

But this follows from the construction of the map $\remb$
and \lemref{vector kitchen}

\end{proof}

Since the roles of $\tg_Q$ and $\tg'_{Q'}$ in $\grD_{G,Q}$ are
essentially symmetric, we have an analog of \propref{catrep}
with $\tg_Q$ replaced by $\tg'_{Q'}$.

\bigskip

Let $K\subset G[[t]]$ be a normal group-subscheme of finite codimension
and let ${\mathfrak k}$ denote its Lie algebra. 
By constructon, the extension
$0\to \CC\to \tg'_{Q'}\to \grg((t))\to 0$ splits canonically
over $\grg[[t]]\subset \grg((t))$. In particular, for every chiral
$\grD_{G,Q}$-module $\M$ supported at the point $x$, we obtain a right
${\mathfrak k}$-action on the corresponding vector space $\M_x$.

\medskip

We let
$\grD_{G,Q}^K\on{-mod}$ denote the full abelian subcategory
of the category of chiral $\grD_{G,Q}$-modules supported at 
$x$, for which the above {\it right} ${\mathfrak k}$-action
integrates to an action of the group-scheme $K$.

Our goal now is to describe the category $\grD_{G,Q}^K\on{-mod}$
in geometric terms. To do that we need to make a digression on
the theory of D-modules on ind-schemes. In what follows we put $Q=0$ 
in order to deal with D-modules (and not with twisted D-modules). 
The generalization to the case of an arbitrary $Q$ is straightforward.

\ssec{D-modules on ind-schemes}

Consider the quotient $Y:=K\backslash G((t))$, as a sheaf 
of sets on the category of quasi-compact schemes with respect
to the faithfully flat topology. 

It is known, cf. \cite{BD1}, 
Theorem 4.5.1, that $Y$ is in fact an ind-scheme of ind-finite type. 
This means that, as a functor, $Y$ is isomorphic to a direct limit 
$\underset{\underset{i\in \NN}{\longrightarrow}}{lim}\,(\underline{Y_i})$, 
where $\underline{Y_i}$ is a functor representable by a scheme $Y_i$ 
and the maps $Y_i\overset{k_{i,j}}\to Y_{j}$ for $j\geq i$ 
are closed embeddings. The ind-finite type property means that the family
$Y_i$ consists of schemes of finite type.

For an ind-scheme $Y$ one defines the category of
$\O^!$-modules as the ind-completion of the direct limit category
$\underset{\longrightarrow}{lim}\,(\O\on{-mod}(Y_i))$, where
for $j\geq i$, the functor $\O\on{-mod}(Y_i)\to \O\on{-mod}(Y_{j})$ 
is of course the direct image under the closed embedding $k_{i,j}$.

\medskip

When $Y$ is of ind-finite type, one can define the category
of {\it right} D-modules on it:

Recall first that if $Y$ is a scheme of finite type the category 
of right D-modules on it makes sense, cf. \cite{BD1}, Sect. 7.10.
Moreover, we have a natural forgetful functor 
$F_Y:\on{D-mod}(Y)\to \O^!\on{-mod}(Y)$. 
(Note that the functor $F_Y$ is left exact but in general
not right exact. If $Y$ is smooth, then it is right exact as
well.) When $Y\overset{k}\to Y'$ is a closed embedding of schemes, there is 
a natural exact functor $k_!:\on{D-mod}(Y)\to \on{D-mod}(Y')$. 
Moreover, we have a natural 
transformation $k_*\circ F_Y\Rightarrow F_{Y'}\circ k_!$. 
\footnote{To distinguish the $\O$- and the D-module direct 
images, we denote the former by $k_*$ and the latter $k_!$.
By * we will always mean the $\O$-module inverse image.}

This allows us to define the category of right D-modules 
on an ind-scheme of finite type. Namely, for 
$Y=\underset{\longrightarrow}{lim}(Y_i)$ as above, we set
$\on{D-mod}(Y)$ to be the ind-completion of the direct limit
of abelian categories
$\underset{\longrightarrow}{lim}\,(\on{D-mod}(Y_i))$.

We have the forgetful functor 
$F_Y:\on{D-mod}(Y)\to \O^!\on{-mod}(Y)$, determined uniquely by the
property that for $\S\in \on{D-mod}(Y_i)$,
$F_Y(\S)=\underset{\underset{j\geq i}{\longrightarrow}}{lim} \,
F_{Y_j}(k_{i,j}{}_!(\S))$. If $Y$ is {\it formally smooth}, cf.
\cite{BD1}, Sect. 7.11.1, the functor $F_Y$ is exact.

\medskip

The main result of this section is the following theorem.
\footnote{We would like to thank A.Beilinson once again
for explaining that \thmref{equivalence} should be
a consequence of the construction of $\grD_{G,Q}$.}

\begin{thm} \label{equivalence}
Let $Y\simeq K\backslash G((t))$ as above.
There is an equivalence of categories
$\on{Sec}:\on{D-mod}(Y)\to\grD_{G,0}^K\on{-mod}$.
\end{thm}

\ssec{Construction of the functor}

Let $\pi$ denote the natural projection $G((t))\to Y$.
For a (right) D-module $\S$ on $Y$ consider the corresponding
$\O^!$-module $F_Y(\S)$. Then the pull-back $\pi^*(F_Y(\S))$ is an 
$\O^!$-module on $G((t))$. Observe that since $G((t))$ is affine,
an $\O^!$-module on it is the same a discrete continuous module over the 
topological algebra $\O_{G((t))}$.

We claim that the vector space $\on{Sec}(\S):=\Gamma(G((t)),\pi^*(F_Y(\S)))$
underlies a chiral $\grD_{G,0}$-module supported at $x\in X$.

First, by construction, $\on{Sec}(\M)$ is a discrete 
$\O_{G((t))}$-module.
Secondly, we endow it with a continuous action of $\grg((t))$ as follows:

Since the projection $\pi$ is right $G$-invariant, 
the right D-module structure
on $\S$ gives rise to the action of $\grg((t))$ on $\pi^*(\S)$ by
derivations of the $\O^!$-module structure, where $\xi\in \grg((t))$
goes to the vector field $-\xi^l$ on $G((t))$.

Therefore, $\on{Sec}(\S)$ indeed corresponds to a chiral 
$\grD_{G,0}$-module, in view of \propref{catrep}.

To complete the construction, it remains to show
that the action of ${\mathfrak k}$ coming from the embedding $\remb$
is indeed integrable.

Obviously, for any $\O^!$-module $\S'$ on $Y$, the pull-back
$\pi^*(\S')$ is $K$-equivariant with respect to the $K$-action 
on $G((t))$ by left translations. 
We will prove the following proposition:

\begin{prop}  \label{actions coincide}
For $\S\in \on{D-mod}(Y)$ the right action of ${\mathfrak k}$ on
$\on{Sec}(\S)$ coincides with the ${\mathfrak k}$-action coming
from the $K$-equivariant structure on $\pi^*(F_Y(\S))$.
\end{prop}

\begin{proof}

Let us view $\on{Sec}$ as a functor from $\on{D-mod}(Y)$
to the category of $\grD_{G,0}$-modules supported at $x$.

The difference of the two actions of ${\mathfrak k}$ for a given
$\S\in \on{D-mod}(Y)$ is a map
$${\mathbf d}:{\mathfrak k}\to \on{End}_{\CC}(\on{Sec}(\S)),$$
which commutes with both the left $\grg((t))$- and the 
$\O_{G((t))}$-actions on $\on{Sec}(\S)$. 

In other words, ${\mathbf d}$ is a map from 
${\mathfrak k}$ to the endomorphism ring of the functor 
$S\mapsto \on{Sec}(\S)$. Our goal is to prove that ${\mathbf d}\equiv 0$. 

Let us first consider the case when $\S$ is the
$\delta$-function D-module $\delta_1$ on $Y$, i.e.
the direct image of $\CC$ under $\on{pt}\to Y$ corresponding to
the coset of $1\in G((t))$. It is easy to see that the corresponding
$\grD_{G,0}$-module identifies with 
$\on{Ind}_{{\mathfrak k}}^{\grg((t))}(\O_K)$ with the natural 
$\grg((t))$- and $\O_{G((t))}$-actions, cf. \propref{catrep}.
In particular, it is generated by a canonical element
$1_{can}\in \on{Ind}_{{\mathfrak k}}^{\grg((t))}(\O_K)$.

Therefore, the map ${\mathbf d}:{\mathfrak k}\to 
\on{End}(\on{Sec}(\S))$ vanishes for $\S=\delta_1$.
Indeed it is enough to show that ${\mathbf d}(1_{can})=0$,
but this follows immediately from the construction.

\medskip

Observe now that 
the action of $G((t))$ by right translations induces endo-functors
of $\on{D-mod}(Y)$ and $\grD_{G,0}$-mod, which commute with the functor 
$\on{Sec}$ in the natural sense. By construction, for every 
$k\in {\mathfrak k}$, ${\mathbf d}(k)\in \on{End}(\on{Sec})$ 
commutes with this $G((t))$-action. 

In particular, since every $\delta_y$ for $y\in Y$,
can be obtained from $\delta_1$ as a $G$-translate,
we conclude that ${\mathbf d}:{\mathfrak k}\to 
\on{End}(\on{Sec}(\delta_y))$ vanishes for all $y\in Y$.

\medskip

The above observations imply what we need:
 
First, every subscheme of $Y$ admits a Zariski-open cover
over which the morphism $\pi$ admits a section. Thus,
let $Y'$ be a locally closed subscheme of $Y$ and $s$ be a map
$Y'\to G((t))$. Then, any D-module $\S$ on $Y'$ is isomorphic as 
an $\O$-module to $s^*(\pi^*(\S))$.

In particular, ${\mathbf d}$ defines a map from 
${\mathfrak k}$ to the endomorphism ring of the 
the forgetful functor $F_{Y'}:\on{D-mod}(Y')\to \O^!\on{-mod}(Y')$, and
it is known that $\on{End}(F_{Y'})\simeq \O_{Y'}$.

However, since ${\mathbf d}$ ``kills'' all the D-modules
of the form $\delta_y$, the map ${\mathbf d}$ vanishes
identically.

\end{proof}

\ssec{Examples}

Let us now describe the
right D-module $\S_{vac}$ on $Y$ corresponding to the vacuum representation 
${\mathbb V}_{G,Q}$ (this was somewhat implicit in the proof of 
\propref{actions coincide} above):

This $\S_{vac}$ will be the direct image of a 
certain right D-module on $K\backslash G[[t]]$ under the closed embedding
$K\backslash G[[t]]\hookrightarrow Y$. 
The corresponding D-module on $K\backslash G[[t]]$ is constructed as follows:

As an $\O$-module it is isomorphic to $\O_{K\backslash G[[t]]}$
and for the left-invariant vector field 
$\xi^l\in \on{Lie}(K\backslash G[[t]])$ and 
$f\in \O_{K\backslash G[[t]]}$ we set:
$f\cdot \xi^l:=-\on{Lie}_{\xi^l}(f)$.
This defines a D-module structure completely, since left-invariant
vector fields and functions generate the ring of differential operators.

Note that in this construction it was crucial that
$K\backslash G[[t]]$ is a group: on an arbitrary variety $Z$, $\O_Z$ 
carries a canonical {\it left} D-module structure but, 
in general, no {\it right} D-module structure.

\medskip

To produce another set of examples, let us recall, cf.  
\cite{BD1}, Sect. 7.11.6, that on any ind-scheme $Y$ of ind-finite 
type, the forgetful functor $F_Y:\on{D-mod}(Y)\to \O^!\on{-mod}(Y)$ 
has an exact left adjoint, which will denote by $I_Y$. (When $Y$ is 
a smooth variety, $I_Y$ is given by 
$\F\mapsto \F\underset{\O_Y}\otimes D_Y$.)

Thus, let us start with an $\O^!$-module $\F$ on $Y$
and consider the corresponding $\grD_{G,0}$-module
$\on{Sec}(I_Y(\F))$. We claim that it can be described as follows:

The pull-back $\pi^*(\F)$ as an $\O^!$-module
on $G((t))$, equivariant with respect to the $K$-action by
left translations. We regard $\Gamma(G((t)),\pi^*(\F))$ as 
a $K$-module and consider the induced $\tg'_{Q'}$-module
$\on{Ind}_{\mathfrak k}^{\tg'_{Q'}}(\pi^*(\F))$.
(Here $Q'=-Q_0$; recall also that the induction 
$\on{Ind}_{\mathfrak k}^{\tg'_{Q'}}$ is understood in the sense
that $1\in\CC\subset \tg'_{Q'}$ acts on the induced module as the 
identity).  

This module carries a compatible $\O_{G((t))}$-action. Hence, as in 
\propref{catrep}, it is naturally a $\grD_{G,0}$-module. 
The right $K$-action on it is integrable by construction and 
we claim that it can be canonically identified with $\on{Sec}(I_Y(\F))$.

Indeed, we have a tautological map of $\tg'_{Q'}$-modules
$\on{Ind}_{\mathfrak k}^{\tg'_{Q'}}(\pi^*(\F))\to 
\on{Sec}(I_Y(\F))$,
which sends $\pi^*(\F)$ identically to 
$\pi^*(\F)\subset \pi^*(I_Y(\F))$. To prove that it is an 
isomorphism we proceed as follows:

Since $Y$ is formally smooth, cf. \cite{BD1}, Theorem 4.5.1
and Proposition 7.11.8, $F_Y(I_Y(\F))$ carries a canonical filtration
such that the associated graded object identifies with
$\F\underset{\O_Y}\otimes \on{Sym}(T_Y)$, where $T_Y$ is 
the tangent sheaf on $Y$.
The $\tg'_{Q'}$-module 
$\on{Ind}_{\mathfrak k}^{\tg'_{Q'}}(\pi^*(\F))$
carries a canonical filtration as well and the above map is
compatible with filtrations. To prove our assertion, it suffices
to observe that the induced map of the graded objects is
\begin{align*}
&\on{gr}(\on{Ind}_{\mathfrak k}^{\tg'_{Q'}}(\pi^*(\F)))\simeq
\Gamma(G((t)),\pi^*(\F))\underset{\CC}\otimes 
\on{Sym}(\grg((t))/{\mathfrak k})\overset{\sim}\to \\
&\Gamma(G((t)),\pi^*(\F\underset{\O_Y}\otimes\on{Sym}(T_Y)))\simeq
\on{gr}(\on{Sec}(I_Y(\F))).
\end{align*}

\ssec{Proof of \thmref{equivalence}}

First, since $Y$ is formally smooth, the functor $F_Y$ is exact and,
hence, the functor $\on{Sec}:\on{D-mod}(Y)\to \grD_{G,0}^K\on{-mod}$
is exact as well. We claim now that it is fully-faithful:

According to \cite{BD1}, Sect. 7.11.8-9, we can think of D-modules
on $Y$, as of $\O^!$-modules endowed with a compatible right action
of the tangent algebroid $T_Y$. Since the natural map
$\grg((t))\otimes \O_Y\to T_Y$ is surjective, the functor $\on{Sec}$
is fully-faithful.

Hence, it remains to show that $\on{Sec}$ is surjective on objects.
However, every object of $\grD_{G,0}^K\on{-mod}$ can be represented
as a quotient of an object of the form
$$\on{Ind}_{\mathfrak k}^{\tg'_{Q'}}(N),$$
where $N$ is a right-$K$-integrable continuous $\O_{G((t))}$-module.

Since we know that all objects of this form are in the image 
of our functor, the theorem follows.

\medskip

\noindent{\it Remark.}
Let us denote by $\M\mapsto\on{Loc}(\M)$ the quasi-inverse functor of
$\on{Sec}$. Explicitely, it can be described as follows:

Let $\M$ be an object of $\grD_{G,0}^K\on{-mod}$. Then
$\M$ can be thought of as a $K$-equivariant $\O^!$-module on 
$G((t))$, and hence, gives rise to an $\O^!$-module on $Y$. 
It acquires a natural (right) action of the Lie-algebroid 
$\grg((t))\otimes \O_Y$ (from {\it minus} the left $\grg((t))$-action on 
$\M$).

Hence, the content of \thmref{equivalence} is that the above 
$\grg((t))\otimes \O_Y$-action factors through the action of the tangent 
algebroid $T_Y$.

\ssec{The category of D-modules on $G((t))$}  \label{finale}

Let now $K_1\subset G[[t]]$ be another normal subgroup,
contained inside $K$. Denote $Y_1:=K_1\backslash G((t))$.
We have the natural exact pull-back functor 
$$\pi_{Y_1,Y}^\circ=\pi_{Y_1,Y}^![-(\on{dim}(K/K'))]:
\on{D-mod}(Y)\to \on{D-mod}(Y'),$$ 
where $\pi_{Y_1,Y}:Y_1\to Y$
is the canonical projection.

Thus, one possible candidate for the category of D-modules 
on $G((t))$ is the ind-completion of the 
direct limit category 
$$\on{D-mod}\,(G((t))):=\on{ind.comp.}\,(
\underset{\underset{K\subset G[[t]]}
\longrightarrow}{lim}(\on{D-mod}(Y))).$$

Note, however, that the functors $F_{Y_1}\circ \pi_{Y_1,Y}^\circ$ and
$\pi_{Y_1,Y}^*\circ F_Y$ from $\on{D-mod}(Y)$ to $\O^!\on{-mod}(Y_1)$
differ by the determinant line: for a D-module $\S$ on $Y$,
$$F_{Y_1}(\pi_{Y_1,Y}^\circ(\S))\simeq \pi_{Y_1,Y}^*(F_Y(\S))
\underset{\CC}\otimes \Lambda^{\on{dim}({\mathfrak k}/{\mathfrak k'})}
({\mathfrak k}/{\mathfrak k'}).$$

Correspondingly, the functors $\on{Sec}:\on{D-mod}(Y)\to 
\grD_{G,0}^K\on{-mod}$ for different subgroups $K$ are compatible 
only up to a twist by the above $1$-dimensional vector space.

Note that the ind-completion of the direct limit category
$\underset{\underset{K\subset G[[t]]}
\longrightarrow}{lim}(\grD_{G,0}^K\on{-mod})$ identifies
with the category of all chiral $\grD_{G,0}$-modules supported at $x$.

Now, \thmref{equivalence} suggests that the above category of
chiral $\grD_{G,0}$-modules is also a reasonable candidate to be
called the category of D-modules on $G((t))$. However,
as we have just seen, there is no canonical equivalence 
(and, propably, no equivalence at all) between
it and the above category $\on{D-mod}\,(G((t)))$.

\end{document}